\documentclass[10pt,a4paper]{article}
\usepackage{amsfonts}
\usepackage{amssymb}
\usepackage{mathrsfs}
\usepackage{amsthm}
\usepackage{setspace}
\usepackage{color}
\usepackage{bm}

\usepackage{fancyhdr}
\newtheorem{theo}{Theorem}[section]
\newtheorem{prop}[theo]{Proposition}
\newtheorem{cor}[theo]{Corollary}
\newtheorem{lemm}[theo]{Lemma}

\newtheorem{con}{Conjecture}

\theoremstyle{definition}

\newenvironment{lproof}{\emph{Proof of Lemma.}}{ \qed \par}

\newcommand{\be}{\begin{eqnarray*}}
\newcommand{\ee}{\end{eqnarray*}}
\newcommand{\beqa}{\begin{eqnarray}}
\newcommand{\eeqa}{\end{eqnarray}}
\newcommand{\ba}{\begin{array}}
\newcommand{\ea}{\end{array}}

\newcommand{\mc}{\mathcal}
\newcommand{\mf}{\mathfrak}

\newcommand{\mbb}{\mathbb}

\newcommand\Dfourcont{\begin{picture}(49,12)\put(24,3){\line(-1,0){18}}%
\put(28,4){\line(2,1){17}}%
\put(27,2){\line(2,-1){18}}%
\put(5,2){\makebox(0,0){$\circ$}}%
\put(26,3){\makebox(0,0){$\circ$}}%
\put(47,12){\makebox(0,0){$\times$}}%
\put(47,-7){\makebox(0,0){$\times$}}%
\end{picture}}

\newcommand\Focto{\begin{picture}(56,12)%
\put(5,3){\line(1,0){16}}%
\put(25,2){\line(1,0){16}}%
\put(25,4){\line(1,0){16}}%
\put(45,3){\line(1,0){16}}%
\put(3,3){\makebox(0,0){$\times$}}%
\put(23,3){\makebox(0,0){$\circ$}}%
\put(43,3){\makebox(0,0){$\circ$}}%
\put(63,3){\makebox(0,0){$\circ$}}%
\put(33,3){\makebox(0,0){$<$}}%
\end{picture}}

\newcommand{\Anum}{\begin{picture}(96,16)%
\put(5,3){\line(1,0){16}}%
\put(25,3){\line(1,0){16}}%
\put(45,3){\line(1,0){6}}%
\put(71,3){\line(-1,0){6}}%
\put(74,3){\line(1,0){17}}%
\put(59,3){\makebox(0,0){\dots}}%
\put(3,3){\makebox(0,0){$\circ$}}%
\put(3,9){\makebox(0,0){$\color{blue}{\scriptstyle{1}}$}}%
\put(23,3){\makebox(0,0){$\circ$}}%
\put(23,9){\makebox(0,0){$\color{blue}{\scriptstyle{1}}$}}%
\put(43,3){\makebox(0,0){$\circ$}}%
\put(43,9){\makebox(0,0){$\color{blue}{\scriptstyle{1}}$}}%
\put(73,3){\makebox(0,0){$\circ$}}%
\put(73,9){\makebox(0,0){$\color{blue}{\scriptstyle{1}}$}}%
\put(93,3){\makebox(0,0){$\circ$}}%
\put(93,9){\makebox(0,0){$\color{blue}{\scriptstyle{1}}$}}%
\end{picture}}

\newcommand{\Bnum}{\begin{picture}(76,16)%
\put(5,3){\line(1,0){16}}%
\put(25,3){\line(1,0){6}}%
\put(51,3){\line(-1,0){6}}%
\put(55,4){\line(1,0){17}}%
\put(54,2){\line(1,0){18}}%
\put(39,3){\makebox(0,0){\dots}}%
\put(63,3){\makebox(0,0){$>$}}%
\put(3,3){\makebox(0,0){$\circ$}}%
\put(3,9){\makebox(0,0){$\scriptstyle{1}$}}%
\put(23,3){\makebox(0,0){$\circ$}}%
\put(23,9){\makebox(0,0){$\color{red}{\scriptstyle{2}}$}}%
\put(53,3){\makebox(0,0){$\circ$}}%
\put(53,9){\makebox(0,0){$\color{red}{\scriptstyle{2}}$}}%
\put(73,3){\makebox(0,0){$\circ$}}%
\put(73,9){\makebox(0,0){$\color{red}{\scriptstyle{2}}$}}%
\end{picture}}

\newcommand{\Cnum}{\begin{picture}(76,16)%
\put(5,3){\line(1,0){16}}%
\put(25,3){\line(1,0){6}}%
\put(51,3){\line(-1,0){6}}%
\put(55,4){\line(1,0){17}}%
\put(54,2){\line(1,0){18}}%
\put(39,3){\makebox(0,0){\dots}}%
\put(63,3){\makebox(0,0){$<$}}%
\put(3,3){\makebox(0,0){$\circ$}}%
\put(3,9){\makebox(0,0){$\color{red}{\scriptstyle{2}}$}}%
\put(23,3){\makebox(0,0){$\circ$}}%
\put(23,9){\makebox(0,0){$\color{red}{\scriptstyle{2}}$}}%
\put(53,3){\makebox(0,0){$\circ$}}%
\put(53,9){\makebox(0,0){$\color{red}{\scriptstyle{2}}$}}%
\put(73,3){\makebox(0,0){$\circ$}}%
\put(73,9){\makebox(0,0){$\scriptstyle{1}$}}%
\end{picture}}

\newcommand{\Dnum}{\begin{picture}(76,17)%
\put(5,4){\line(1,0){16}}%
\put(25,4){\line(1,0){6}}%
\put(51,4){\line(-1,0){6}}%
\put(55,5){\line(5,1){17}}%
\put(55,3){\line(5,-1){18}}%
\put(39,4){\makebox(0,0){\dots}}%
\put(3,4){\makebox(0,0){$\circ$}}%
\put(3,10){\makebox(0,0){$\color{blue}{\scriptstyle{1}}$}}%
\put(23,4){\makebox(0,0){$\circ$}}%
\put(23,10){\makebox(0,0){$\color{red}{\scriptstyle{2}}$}}%
\put(53,4){\makebox(0,0){$\circ$}}%
\put(53,10){\makebox(0,0){$\color{red}{\scriptstyle{2}}$}}%
\put(74,8){\makebox(0,0){$\circ$}}%
\put(78,9){\makebox(0,0){$\color{blue}{\scriptstyle{1}}$}}%
\put(74,0){\makebox(0,0){$\circ$}}%
\put(78,0){\makebox(0,0){$\color{blue}{\scriptstyle{1}}$}}%
\end{picture}}

\newcommand{\Gnum}{\begin{picture}(26,16)%
\put(5,3){\line(1,0){16}}%
\put(4,4){\line(1,0){18}}%
\put(4,2){\line(1,0){18}}%
\put(13,3){\makebox(0,0){$<$}}%
\put(3,3){\makebox(0,0){$\circ$}}%
\put(3,9){\makebox(0,0){$\scriptstyle{3}$}}%
\put(23,3){\makebox(0,0){$\circ$}}%
\put(23,9){\makebox(0,0){$\color{red}{\scriptstyle{2}}$}}%
\end{picture}}

\newcommand{\Fnum}{\begin{picture}(56,16)%
\put(5,3){\line(1,0){16}}%
\put(25,2){\line(1,0){16}}%
\put(25,4){\line(1,0){16}}%
\put(45,3){\line(1,0){16}}%
\put(3,3){\makebox(0,0){$\circ$}}%
\put(3,9){\makebox(0,0){$\color{red}{\scriptstyle{2}}$}}%
\put(23,3){\makebox(0,0){$\circ$}}%
\put(23,9){\makebox(0,0){$\scriptstyle{4}$}}%
\put(43,3){\makebox(0,0){$\circ$}}%
\put(43,9){\makebox(0,0){$\scriptstyle{3}$}}%
\put(63,3){\makebox(0,0){$\circ$}}%
\put(63,9){\makebox(0,0){$\color{red}{\scriptstyle{2}}$}}%
\put(33,3){\makebox(0,0){$<$}}%
\end{picture}}

\newcommand{\Esixnum}{\begin{picture}(86,34)%
\put(5,3){\line(1,0){16}}%
\put(25,3){\line(1,0){16}}%
\put(45,3){\line(1,0){16}}%
\put(43,5){\line(0,1){16}}%
\put(65,3){\line(1,0){16}}%
\put(3,3){\makebox(0,0){$\circ$}}%
\put(3,9){\makebox(0,0){$\color{blue}{\scriptstyle{1}}$}}%
\put(23,3){\makebox(0,0){$\circ$}}%
\put(23,9){\makebox(0,0){$\color{red}{\scriptstyle{2}}$}}%
\put(43,3){\makebox(0,0){$\circ$}}%
\put(46,9){\makebox(0,0){$\scriptstyle{3}$}}%
\put(43,23){\makebox(0,0){$\circ$}}%
\put(43,29){\makebox(0,0){$\color{red}{\scriptstyle{2}}$}}%
\put(63,3){\makebox(0,0){$\circ$}}%
\put(63,9){\makebox(0,0){$\color{red}{\scriptstyle{2}}$}}%
\put(83,3){\makebox(0,0){$\circ$}}%
\put(83,9){\makebox(0,0){$\color{blue}{\scriptstyle{1}}$}}%
\end{picture}}

\newcommand{\Esevennum}{\begin{picture}(106,34)%
\put(5,3){\line(1,0){16}}%
\put(25,3){\line(1,0){16}}%
\put(45,3){\line(1,0){16}}%
\put(43,5){\line(0,1){16}}%
\put(65,3){\line(1,0){16}}%
\put(85,3){\line(1,0){16}}%
\put(3,3){\makebox(0,0){$\circ$}}%
\put(3,9){\makebox(0,0){$\color{red}{\scriptstyle{2}}$}}%
\put(23,3){\makebox(0,0){$\circ$}}%
\put(23,9){\makebox(0,0){$\scriptstyle{3}$}}%
\put(43,3){\makebox(0,0){$\circ$}}%
\put(46,9){\makebox(0,0){$\scriptstyle{4}$}}%
\put(43,23){\makebox(0,0){$\circ$}}%
\put(43,29){\makebox(0,0){$\color{red}{\scriptstyle{2}}$}}%
\put(63,3){\makebox(0,0){$\circ$}}%
\put(63,9){\makebox(0,0){$\scriptstyle{3}$}}%
\put(83,3){\makebox(0,0){$\circ$}}%
\put(83,9){\makebox(0,0){$\color{red}{\scriptstyle{2}}$}}%
\put(103,3){\makebox(0,0){$\circ$}}%
\put(103,9){\makebox(0,0){$\scriptstyle{1}$}}%
\end{picture}}

\newcommand{\Eeightnum}{\begin{picture}(126,34)%
\put(5,3){\line(1,0){16}}%
\put(25,3){\line(1,0){16}}%
\put(45,3){\line(1,0){16}}%
\put(43,5){\line(0,1){16}}%
\put(65,3){\line(1,0){16}}%
\put(85,3){\line(1,0){16}}%
\put(105,3){\line(1,0){16}}%
\put(3,3){\makebox(0,0){$\circ$}}%
\put(3,9){\makebox(0,0){$\color{red}{\scriptstyle{2}}$}}%
\put(23,3){\makebox(0,0){$\circ$}}%
\put(23,9){\makebox(0,0){$\scriptstyle{4}$}}%
\put(43,3){\makebox(0,0){$\circ$}}%
\put(46,9){\makebox(0,0){$\scriptstyle{6}$}}%
\put(43,23){\makebox(0,0){$\circ$}}%
\put(43,29){\makebox(0,0){$\scriptstyle{3}$}}%
\put(63,3){\makebox(0,0){$\circ$}}%
\put(63,9){\makebox(0,0){$\scriptstyle{5}$}}%
\put(83,3){\makebox(0,0){$\circ$}}%
\put(83,9){\makebox(0,0){$\scriptstyle{4}$}}%
\put(103,3){\makebox(0,0){$\circ$}}%
\put(103,9){\makebox(0,0){$\scriptstyle{3}$}}%
\put(123,3){\makebox(0,0){$\circ$}}%
\put(123,9){\makebox(0,0){$\color{red}{\scriptstyle{2}}$}}%
\end{picture}}

\newcommand{\Cquater}{\begin{picture}(76,12)%
\put(5,3){\line(1,0){16}}%
\put(25,3){\line(1,0){6}}%
\put(51,3){\line(-1,0){6}}%
\put(55,4){\line(1,0){17}}%
\put(54,2){\line(1,0){18}}%
\put(39,3){\makebox(0,0){\dots}}%
\put(63,3){\makebox(0,0){$<$}}%
\put(3,3){\makebox(0,0){$\circ$}}%
\put(23,3){\makebox(0,0){$\times$}}%
\put(53,3){\makebox(0,0){$\circ$}}%
\put(73,3){\makebox(0,0){$\circ$}}%
\end{picture}}

\newcommand{\DfreeCR}{\begin{picture}(76,12)%
\put(5,3){\line(1,0){16}}%
\put(25,3){\line(1,0){6}}%
\put(51,3){\line(-1,0){6}}%
\put(55,4){\line(2,1){17}}%
\put(54,2){\line(2,-1){18}}%
\put(39,3){\makebox(0,0){\dots}}%
\put(3,3){\makebox(0,0){$\circ$}}\put(23,3){\makebox(0,0){$\circ$}}%
\put(53,3){\makebox(0,0){$\circ$}}%
\put(74,12){\makebox(0,0){$\times$}}%
\put(73,-7){\makebox(0,0){$\times$}}%
\put(74,2){\makebox(0,0){$\updownarrow$}}%
\end{picture}}

\newcommand{\DfreeLan}{\begin{picture}(76,12)%
\put(5,3){\line(1,0){16}}%
\put(25,3){\line(1,0){6}}%
\put(51,3){\line(-1,0){6}}%
\put(55,4){\line(2,1){17}}%
\put(54,2){\line(2,-1){18}}%
\put(39,3){\makebox(0,0){\dots}}%
\put(3,3){\makebox(0,0){$\circ$}}%
\put(23,3){\makebox(0,0){$\circ$}}%
\put(53,3){\makebox(0,0){$\circ$}}%
\put(74,12){\makebox(0,0){$\times$}}%
\put(73,-7){\makebox(0,0){$\times$}}%
\end{picture}}

\begin{document}

\title{Non-regular $|2|$-graded geometries II: classifying geometries, and generic six-in-nine distributions}
\author{Stuart Armstrong, \\ St Cross College, \\ Oxford, OX1 3LZ, UK}
\date{2009}
\maketitle

\begin{abstract}
Complementing the previous paper in the series, this paper classifies $|2|$-graded parabolic geometries, listing their important properties: the group $G_0$, the graded tangent bundle $gr(T)$ and its algebra\"ic bracket, the relevant cohomology spaces and the standard Tractor bundle $\mc{T}$. Several of these geometries are then explored in more detail, and the paper ends with a case study that that partially solves the equivalence problem for generic six distributions on nine dimensional manifolds.
\end{abstract}

\section{Introduction}
Two groups $G$ and $P$ form a $|k|$-graded parabolic pair if $G$ is semisimple, and if the Lie algebra $\mf{g}$ of $G$ admits a $|k|$-grading
\be
\mf{g} &=& \mf{g}_{-k} \ \oplus \ \ldots \ \oplus \ \mf{g}_{-1} \ \oplus \ \mf{g}_{0} \ \oplus \ \mf{g}_{1} \ \oplus \ \ldots \ \oplus \mf{g}_k,
\ee
such that $[\mf{g}_i, \mf{g}_j] \subset \mf{g}_{i+j}$ and such that the Lie algebra $\mf{p}$ of $P$ is given by
\be
\mf{p} = \sum_{i \geq 0} \mf{g}_i.
\ee
The reductive part of $P$ is the group $G_0$, which has Lie algebra $\mf{g}_0$. The nilpotent part has Lie algebra $\mf{g}^+ = \sum_{i \geq 1} \mf{g}_i$.

For $(G,P)$ a $|2|$-graded parabolic pair, the previous paper \cite{me2grad1} defined partially regular $|2|$-graded geometries, namely a manifold $M$, with a distribution $H$ contained in the tangent bundle $T$ and an reduction of the structure bundle of $gr(T) = T/H \oplus H$ to $\mc{P}$, a $P$-principal bundle, such that
\be
gr(T) = \mc{P} \times_P (\mf{g}/\mf{p}).
\ee

The $P$ structure defines an algebra\"ic bracket $\mc{K} \in \Gamma(\wedge^2 H^* \otimes (T/H))$, and a corresponding $\partial^*$ bundle operator
\be
\partial^*: \ \wedge^2 H^* \otimes (T/H) \ \ \to \ \ (H^* \otimes H) \oplus ((T/H)^* \otimes (T/H)).
\ee
There is another natural section of $\wedge^2 H^* \otimes (T/H)$, the Levi bracket $\mc{L}$, defined for sections $X$ and $Y$ of $H$ as
\be
\mc{L}(X,Y) = [X,Y] / H.
\ee
The geometry is then said to be partially regular if
\be
\partial^* (\mc{K} - \mc{L}) = 0.
\ee

The main result was
\begin{theo}
If the second cohomology space $H^1(\mf{g}^+, \mf{g})$ vanishes in strictly positive homogeneity, there exists a unique normal Cartan connection encoding the geometry.
\end{theo}

Furthermore, if the only data is a generic holomorphic distribution $H$ and a $|2|$-graded complex parabolic pair $(G^{\mbb{C}},P^{\mbb{C}})$ such that $dim_{\mbb{C}} \ \mf{g}_{-1} = rank_{\mbb{C}} \ H$ and $dim_{\mbb{C}} \ \mf{g}_{-2} = rank_{\mbb{C}} \ T/H$, then
\begin{theo} \label{second:theo}
Almost everywhere on $M$, there exists a principal bundle $\mc{P}^{\mbb{C}}$ for $gr(T)$ with structure group $P^{\mbb{C}}$, such that $(M,H,\mc{P}^{\mbb{C}})$ is a partially regular $|2|$-graded geometry.

If $H^1(\mf{g}^+, \mf{g})$ vanishes in homogeneity zero, there is a finite choice of this bundle $\mc{P}^{\mbb{C}}$, up to isomorphism.
\end{theo}
If $\kappa$ is the curvature of a regular normal Cartan connection, then its lowest homogeneity component is $\partial$-closed, hence a section of
\be
\mc{P} \times_P H^2(\mf{g}^+,\mf{g}).
\ee
This property extends to partially regular normal $|2|$-graded Cartan connections.

Hence the important formal elements of a $|2|$-graded geometry connection are:
\begin{enumerate}
\item the group $G_0$,
\item the bundle $gr(T)$,
\item the cohomology spaces $H^1(\mf{g}^+, \mf{g})$ (in non-negative homogeneities) and $H^2(\mf{g}^+, \mf{g})$, and
\item the algebra\"ic bracket $\mc{K}$.
\end{enumerate}
Also of interest for Tractor calculus (see \cite{TCPG} and \cite{TBIPG}) is the standard Tractor bundle $\mc{T}$.

This paper will list these for all simple $|2|$-graded complex pairs $(G,P)$. Some of the real forms of these geometries are of particular interest: non-integral quaternionic- and octonionic-contact geometries, free-Lagrangians and free-CR geometries, and conformal-spin geometries. These will be analysed seperately.

Finally, the `almost everywhere' portion of Theorem \ref{second:theo} is dependent only the value of the Levi bracket at a point. Though paper \cite{me2grad1} establishes that almost all Levi brackets are suitable, it gives no criteria for deciding if a given bracket works. This paper concludes by presenting the full solution to this problem for distributions of rank six in a nine-dimensional manifold, both in the real and the holomorphic categories, and presents a tractable sufficient condition to for analysing whethere a given $H$ can be partially regularised. It does this by first demonstrating that the orbit classification of \cite{metabelian} extends to the real category, then by identifying the families of orbits that a finite choice of allow partial regularisations.

This defines precisely what is meant by a \emph{generic} distribution of rank six-in-nine and -- in theory -- solves the equivalence problem for such distributions, and establishes (\cite{capkatauto} and \cite{capauto}) that their automorphism group is a Lie group of dimension bounded by $28$, the dimension of the group $G = Spin(4,4)$.

\subsection*{Acknowledgements}
It gives me great pleasure to acknowledge the financial support project P19500-N13 of the ``Fonds zur F\"orderung der wissenschaftlichen Forschung (FWF)'', as well as the help and comments of Andreas {\v{C}}ap and Jan Slov{\'a}k, and the support of the math department of the university of Vienna.

\section{Simple $|2|$-graded geometries}

Initially, we will assume that the group $G$ is simple. Semisimple examples can then be built up by combining these simple examples with other $|1|$- or $|2|$-graded geometries; the process will be described in section \ref{semi:simp}.

The complexified versions of these geometries are constructed by taking a Dynkin diagram and crossing a certain number of roots \cite{capslo}. Let $\tau$ be the maximal torus, $\Phi$ the set of all simple roots, $\Phi_+$ the set of all simple positive roots and $\Phi_- = -\Phi_+$ the set of negative simple roots. Now let $R_- \subset \Phi_-$ consist of all negative simple roots that are crossed on the Dynkin diagram, and similarly define $R_+ = - R_- \subset \Phi_+$.

Then $\mf{p}$ may be defined as the direct sum of $\tau$ with all root spaces for the roots of type
\be
\sum_{\alpha \in \Phi, \alpha \notin R_-} n_{\alpha} \alpha,
\ee
for natural numbers $n_{\alpha}$. Note that when all roots are crossed, $R_- = \Phi_-$ and $\mf{p}$ is the Borel subalgebra.

The group $P$ can then be defined as the subgroup of $G$ that acts as automorphisms on $\mf{p}$. The algebra $\mf{g}_0$ is spanned $\tau$ and by root spaces for roots of type
\be
\sum_{{\alpha \in \Phi, \alpha \notin (R_+ \cup R_-)}} n_{\alpha} \alpha.
\ee

The highest weight of the adjoint representation of $G$ is equal to a sum $\sum_{\alpha \in \Phi_+} m_{\alpha} \alpha$ for certain natural numbers $m_{\alpha} \neq 0$. Then the geometry $(G,P)$ is $|k|$-graded, where
\be
k = \sum_{\alpha \in R_+} m_{\alpha}.
\ee
These $m_{\alpha}$ are known; table \ref{table:numbers} presents them for each simple Dynkin diagram.

\begin{table}[htbp]
\begin{center}
\begin{tabular}{|c|c|}
\hline $A_n$ & \Anum \\
\hline $B_n$ & \Bnum \\
\hline $C_n$ & \Cnum \\
\hline $D_n$ & \Dnum \\
\hline $G_2$ & \Gnum \\
\hline $F_4$ & \Fnum \\
\hline $E_6$ & \Esixnum \\
\hline $E_7$ & \Esevennum \\
\hline $E_8$ & \Eeightnum \\
\hline
\end{tabular}
\end{center}
\caption{$m_{\alpha}$ for each root}
\label{table:numbers}
\end{table}

Simple $|2|$-graded geometries are thus generated by crossing a single double root (those where $m_{\alpha} = 2$), or a pair of single roots ($m_{\alpha} = 1$). Section \ref{class:results} will present all such simple geometries, defining the main characteristics: the bundle $T_{-1}$, the bundle $T_{-2} = T/T_{-1}$, the group $G_0$, the form of the algebra\"ic bracket $\mc{K}$ between $\wedge^2 T_{-1}$ and $T_{-2}$, the standard Tractor bundle $\mc{T}$, and the relevant cohomology spaces.

For simplicity, we will only deal with the holomorphic forms of parabolic geometries, with $G$ and $P$ both complex. The various real forms and their corresponding geometries can be deduced from this. In general, the bundle and bracket information is easy to deduce, from simply considering the dimensions of $\mf{g}$ and $\mf{g}^+$, the dimensions of possible $G_0$-modules, the fact that $G_0$ must act faithfully on $\mf{g}_{-1}$, and the possible inclusions of  $\mf{g}_{-2}$ into $\wedge^2 \mf{g}_{-1}$.

This is not enough for $D_{m}$ family, where the above methods do not fully characterise the geometry for $\mf{g}_0 = \mbb{C}^2 \oplus sl(m-1, \mbb{C})$; there are two possibilities. In fact, both of them actually occur, giving the free-Lagrangian geometry and the conformal spin structures respectively -- see sections \ref{free:la} and \ref{con:spin}.

This is also insufficient for $E_6^{\mbb{C}}$, in the case where $\mf{g}_0 = \mbb{C} \oplus \mf{sl}(5,\mbb{C}) \oplus \mf{sl}(2,\mbb{C})$. There are two possibilities for $\mf{g}_-$. Namely, both $(\odot^2 \mbb{C}^5) \oplus \big(\mbb{C}^2 \otimes \mbb{C}^5 \big)$ and $(\mbb{C}^5)^*  \oplus \big( \mbb{C}^2 \otimes ( \wedge^2 \mbb{C}^5) \big)$ are of the right dimensions and define a palausible Lie bracket. The first possibility is excluded, however, by careful analysis of fundamental weight multiplicities within the roots of $E_6^{\mbb{C}}$.

The cohomology is computed using Kostant's proof of the Bott-Borel-Weyl theorem, see \cite{Kostant}. There is one special class of geometries that can be set aside now: Lie contact geometries (see \cite{liecontact2} and \cite{liecontact}). These are geometries where the crossed root(s) on the Dynkin diagram correspond to the root(s) for the adjoint representation (first and last root for the $A_m$ family, second root for the $B_m$ and $D_m$ families, first root for the $C_m$ families, and see section \ref{excep:geo} for the exceptional geometries). They all have $T_{-2}$ being rank one, making $T_{-1}$ into a contact distribution.

In homogeneity zero, contact geometries all have maximally non-vanishing $H^1(\mf{g}^+,\mf{g})$ and vanishing $H^2(\mf{g}^+,\mf{g})$. The first result simply encodes the fact that every contact distribution is locally isomorphic to every other: hence the extra structure on the geometry cannot be deduced from the distribution. There may be a certain analogue to non-regular geometry for CR structures (see \cite{meCR}), but these are very different objects to the geometries we are considering here.

The cohomology results for simple $G$ can be summarised in the following theorems:
\begin{theo} \label{first:cohmo}
For $(\mf{g},\mf{p})$ a simple, $|2|$-graded parabolic pair, the first cohomology space $H^1(\mf{g}^+,\mf{g})$ always vanishes in strictly positive homogeneities. It will vanish in homogeneity zero unless $(\mf{g},\mf{p})$ is a contact geometry, or is from the $A_m$ family with either the first or last root crossed.
\end{theo}
This is a simple consequence of \cite{Kostant} (see also \cite{capslo}). Consequently, almost all $|2|$-graded geometries are determined entirely by the distribution $T_{-1}$.

If a non-regular normal geometry is to exist, there must be a piece of harmonic curvature in homogeneity zero (see \cite{me2grad1}). Consequently, there must a homogeneity zero piece of the second cohomology. In general, this is the case:
\begin{theo} \label{may:be:nonreg}
For $(\mf{g},\mf{p})$ a simple, $|2|$-graded parabolic pair, there is a homogeneity zero piece of $H^2(\mf{g}^+,\mf{g})$ unless
\begin{itemize}
\item $(\mf{g},\mf{p})$ is a contact geometry,
\item $\mf{g} = A_3$, with the first and second (or second and third) roots crossed,
\item $\mf{g} = A_4$, with the second and third roots crossed.
\item $\mf{g} = B_m$, with the last root crossed,
\item $\mf{g} = C_3$, with the second root crossed.
\end{itemize}
\end{theo}
So almost all $|2|$-graded geometries can admit non-regular versions. If, however, the second cohomology lies entirely in homogeneity zero, we have a stronger result: if the geometry is flat, then it must be non-regular. This is also the case for all most all $|2|$-graded geometries:
\begin{theo} \label{must:be:nonreg}
For $(\mf{g},\mf{p})$ a simple, $|2|$-graded parabolic pair, $H^2(\mf{g}^+,\mf{g})$ is concentrated entirely in homogeneity zero unless
\begin{itemize}
\item $(\mf{g},\mf{p})$ is one of the cases listed in Theorem \ref{first:cohmo} or Theorem \ref{may:be:nonreg},
\item $\mf{g} = A_m$ with the second or second-to-last root crossed,
\item $\mf{g} = A_m$ with two adjacent roots crossed,
\item $\mf{g} = B_m $ with the third root crossed,
\item $\mf{g} = C_m $ with the second or second-to-last root crossed,
\item $(\mf{g},\mf{p})$ is a conformal-spin geometry \ref{con:spin},
\item $\mf{g} = D_m$ with the third root crossed.
\end{itemize}
\end{theo}
Thus in general, all non-flat $|2|$-graded geometries must be non-regular -- this includes, for instance, all the non-contact exceptional geometries.

\section{Semisimple $|2|$-graded geometries} \label{semi:simp}
Semisimple $|2|$-graded geometries are modelled on parabolics $G/P$ where
\be
\mf{g} &=& \oplus_i \ \mf{g}^i, \\
\mf{p} &=& \oplus_i \ \mf{p}^i,
\ee
such that $(\mf{g}^i, \mf{p}^i)$ is a $|1|$- or $|2|$-graded parabolic pair, and at least one of these pairs is $|2|$-graded.

Let us first concentrate on the situations where we have two summands (not necessarily simple), $\mf{g} = \mf{g}^1 \oplus \mf{g}^2$ and $\mf{p} = \mf{p}^1 \oplus \mf{p}^2$.

Then by \cite{capslo}'s treatment of Kostant's version of the Bott-Borel-Weil theorem, we can affirm that
\beqa \label{hom:decom}
H^n (\mf{g}^1, \mf{g}) &\cong \oplus_{i+j = n} & \big( H^i ((\mf{g}^1)^+, \mf{g}^1) \otimes H^j ((\mf{g}^2)^+, \mbb{C}) \\
\nonumber && \oplus  \ H^i ((\mf{g}^1)^+, \mbb{C}) \otimes  H^i ((\mf{g}^2)^+, \mf{g}^2) \big)
\eeqa
Standard cohomologies are also calculated in \cite{capslo}, namely the fact that $H^0((\mf{g}^i)^+, \mbb{C}) = \mbb{C}$ and $H^0((\mf{g}^i)^+, \mf{g}^i) = \mf{g}_{-k}$, where $k$ is the grading of $(\mf{g}_i, \mf{p}_i)$. Furthermore, $H^1((\mf{g}^i)^+,\mbb{C}) = \mf{g}_1^i$. Hence we can tackle the first cohomology in the semisimple case:

\begin{prop} \label{hom:decom:prop}
If all of the summands $(\mf{g}^i, \mf{p}^i)$ are $|2|$-graded,
\be
H^1(\mf{g}^1,\mf{g}) = \oplus_i \left( H^1((\mf{g}^i)^+,\mf{g}^i) \right).
\ee
\end{prop}
\begin{proof}
Assume $\mf{g} = \mf{g}^1 \oplus \mf{g}^2$, with both $(\mf{g}^1,\mf{p}^1)$ and $(\mf{g}^2,\mf{p}^2)$ being $|2|$-graded. Then equation (\ref{hom:decom}) implies there are four summands in $H^1(\mf{g}^1,\mf{g})$; two of them are of type
\be
H^1 ((\mf{g}^i)^+, \mf{g}^i) \otimes H^0 ((\mf{g}^j)^+, \mbb{C}) = H^1 ((\mf{g}^i)^+, \mf{g}^i) \otimes \mbb{C} = H^1 ((\mf{g}^i)^+, \mf{g}^i).
\ee
The other two are of type
\be
H^0 ((\mf{g}^i)^+, \mf{g}^i) \otimes H^1 ((\mf{g}^j)^+, \mbb{C}) = \mf{g}_{-k}^i \otimes \mf{g}_{1}^j.
\ee
Since both of these must be of homogeneity $-k+1 = -1$, the equality in homogeneity zero is established.
\end{proof}
We now need to look at the second cohomology. The results aren't strong here, since we have less control over summands of the type $H^2 ((\mf{g}^j)^+, \mbb{C})$. All that can be said with confidence is:
\begin{prop}
\be
H^2(\mf{g}^+,\mf{g}) \supset \oplus_i \left( H^2((\mf{g}^i)^+,\mf{g}^i) \right).
\ee
\end{prop}

Going back to Proposition \ref{hom:decom:prop}, we see that this may give us a way of classifying generic holomorphic distributions of the right rank and co-rank via the methods in \cite{me2grad1}: let $H$ is a generic distribution of rank $r$ and co-rank $s$. Then assume we can construct a parabolic pair
\be
(\mf{g} = \oplus_i \ \mf{g}^i, \mf{p} = \oplus_i \ \mf{p}^i),
\ee
where each of the $(\mf{g}^i, \mf{p}^i)$ is complex $|2|$-graded parabolic pair with vanishing first cohomology in non-negative homogeneities, and with
\be
\textrm{dim}_{\mbb{C}} \ \mf{g}_{-1} &=& r \\
\textrm{dim}_{\mbb{C}} \ \mf{g}_{-2} &=& s.
\ee
Then the corresponding $(\mf{g},\mf{p})$ geometry will depend only on the distribution $T_{-1}$, and, once we have chosen a uniqueness procedure, a generic distribution $T_{-1}$ will be classified almost everywhere by a unique normal Tractor connection. The possible values of $r$ and $s$ can be determined by looking first at the simple examples with the required cohomological condition. These are known from Theorem \ref{may:be:nonreg}, and are listed in table \ref{table:distributons}; in this table, $t_n$ is the $n$-th triangular number, $t_n = n(n-1)/2$.

\begin{table}[htbp]
\begin{center}
\begin{tabular}{|c|c||c|c|}
\hline Co-rank of $T_{-1}$ & Rank of $T_{-1}$ & Co-rank of $T_{-1}$ & Rank of $T_{-1}$  \\
\hline \hline
\hline $\boldsymbol{p \times q} \ \ \ p,q\geq 2$ & $\boldsymbol{l \times (p+q)} \ \ \ l \geq 1$ & $\boldsymbol{5}$ & $\boldsymbol{20}$ \\
\hline $\boldsymbol{t_n} \ \ \ n \geq 3$ & $ \boldsymbol{l\times n} \ \ \ l \geq 1$ & $\boldsymbol{10}$ & $\boldsymbol{32}$ \\
\hline $\boldsymbol{t_{n+1}} \ \ \ n \geq 2$ & $ \boldsymbol{2l\times n} \ \ \ l \geq 1$ & $\boldsymbol{7}$ & $\boldsymbol{35}$\\
\hline $\boldsymbol{7}$ & $\boldsymbol{8}$ & $\boldsymbol{14}$ & $\boldsymbol{63}$ \\
\hline $\boldsymbol{8}$ & $\boldsymbol{16}$ & & \\
\hline
\end{tabular}
\end{center}
\caption{Dimensions and co-dimensions of simple $|2|$-graded parabolic geometries with $H^1(\mf{g}^+, \mf{g}) = 0$ in non-negative homogeneities}
\label{table:distributons}
\end{table}

Since $r$ may be arbitrarily large even for a small $s$, it is easier to initially classify these geometries by the co-rank. The co-rank can never be one or two, and the list of possible ranks for the first few co-ranks are given in table \ref{table:ranks}.
\begin{table}[htbp]
\begin{center}
\begin{tabular}{|c||c|c|}
\hline Co-rank of $T_{-1}$ & Possible rank of $T_{-1}$, general rule & Exceptions\\
\hline
\hline 3 & $3q$ or $4q$ for $q\geq 1$&\\
\hline 4 & $4q$ for $q\geq 1$&\\
\hline 5 & $20$ & \\
\hline 6 & $r \geq 4$ &\\
\hline 7 & $r \geq 7$ & $r \neq 9$ \\
\hline 8 & $6q$ or $(20 + 3q)$ for $q \geq 1$, $4p$ for $p \geq 2$ &\\
\hline 9 & $r \geq 6$ &\\
\hline 10 & $r \geq 5$& $r \neq 6$\\
\hline 11 & $ r \geq 9$ & $r \neq 13$\\
\hline 12 & $r \geq 7$ &\\
\hline 13 & $r \geq 8$ & $r \neq 9$\\
\hline 14 & $r \geq 9$ &\\
\hline 15 & $r \geq 6$ & $r \neq 7,9$\\
\hline 16 & $r \geq 8$ & \\
\hline 17 & $r \geq 12$ &\\
\hline 18 & $r \geq 9$ & \\
\hline 19 & $r \geq 10$ &\\
\hline 20 & $r \geq 9$ & $r \neq 11$\\
\hline 21 & $r \geq 7$ & $ r \neq 8,9$\\
\hline 22 & $r \geq 12$ &\\
\hline 23 & $r \geq 12$ &\\
\hline
\end{tabular}
\end{center}
\caption{Dimensions and co-dimensions of semisimple $|2|$-graded parabolic geometries with $H^1(\mf{g}^+, \mf{g}) = 0$ in non-negative homogeneities}
\label{table:ranks}
\end{table}
Since $[H,H] = T$, we must have $s < r(r-1)/2$. With this in mind, and looking at the values on the table, the tentative conjecture suggests itself.
\begin{con}
The minimal possible rank of $H$ grows as the square root of twice the co-rank; exceptions will remain very rare and may vanish entirely.
\end{con}

\subsection{Classification of simple $|2|$-grade geometries} \label{class:results}
\subsubsection{Linear geometries}
Here $G = SL(m,\mbb{C})$, and two roots are crossed, giving $G_0 = (\mbb{C}^*)^2 \times SL(a,\mbb{C}) \times SL(b,\mbb{C}) \times SL(c,\mbb{C})$, where $a,b$ and $c$ are strictly positive integers with $a+b+c =m$.

The geometry is given by three vector bundles $H_1$, $H_2$ and $K$ of rank $a$, $c$ and $b$, such that
\be
T_{-1} &=& (K \otimes H_1) \oplus (K^* \otimes H_2) \\
T_{-2} &=& H_1 \otimes H_2 \\
\mc{T} &=& H_1^* \oplus K \oplus H_2.
\ee
The algebra\"ic bracket is given by the tensor product on $H_1$ and $H_2$, and the contraction of $K$ with $K^*$.

These geometries have the most complicated cohomological behaviour. Generically, the second cohomology consists of five components, all of homogeneity zero. The following rules then apply:
\begin{enumerate}
\item If one of the crossed roots is first or last, there is one less cohomology component.
\item If the two crossed roots are neighbours, there is one less cohomology component.
\item If one of the crossed roots is first or last, the total homogeneity of all the components goes up by two.
\item If one of the crossed roots is second or second-to-last, the total homogeneity of all the components goes up by one.
\item If the two crossed roots are neighbours, the total homogeneity of all the components goes up by two. \label{end:list} 
\end{enumerate}

These rules are cumulative; for instance, consider $A_2$ with both roots crossed. Then we lose one cohomology component each for having the first and the last root crossed, and one more for having neighbours crossed, so have $5-1-1-1 = 2$ cohomology components in total. As for total homogeneity, we have $+2$ for neighbouring crossed roots, $+2$ for both the first root and the last root and $+1$ for both the second and second-to-last root, for a total of $2 + 2 \times 2+ 2 \times 1 = 8$. This is correct; in fact, the two components are both of homogeneity four.

In order to fully establish the number of pieces and their homogeneities, the next rule is generally sufficient:
\begin{enumerate}
\setcounter{enumi}{5}
\item Whenever the total homogeneity is less than or equal to the number of cohomology components, all components are of homogeneity one or zero.
\end{enumerate}

The cases not covered by the above rule are:
\begin{itemize}
\item $A_2$ with both roots crossed, where the homogeneities are $4$ and $4$.
\item $A_m$, $m>2$ with the first and last roots crossed, where the homogeneities are $2$, $1$ and $1$.
\item $A_3$ with the first and second (or second and last) roots crossed, where the homogeneities are $3$, $2$ and $1$.
\item $A_m$, $m>3$ with the first two (or last two) roots crossed, where the homogeneities are $3$, $2$ and $0$.
\end{itemize}

\subsubsection{Orthogonal geometries}
Here $G=SO(2m+1,\mbb{C})$ or $SO(2m,\mbb{C})$. In the first case, for the $B_{m}$ family, the $|2|$-graded geometries are given by crossing any single root apart from the first one. This gives $G_0 = \mbb{C}^* \times SL(a,\mbb{C}) \times SO(2b+1,\mbb{C})$ where $a+b = m$.

Most of the $D_{m}$ family are of the same type, given by crossing any single root apart from the first one or the last two. This gives $G_0 = \mbb{C}^* \times SL(a,\mbb{C}) \times SO(2b,\mbb{C})$ where $a+b = m$.

In these cases, the geometry is given by two vector bundles $H$ and $K$ of rank $a$ and $2b$ (or $2b+1$), with a metric $g$ on $K$. They have the properties that
\be
T_{-1} &=& K \otimes H \\
T / T_{-1} &=& H \wedge H \\
\mc{T} &=& H^* \oplus K \oplus H.
\ee
The algebra\"ic bracket is given by the wedge on $H$, and the contraction of $K$ with itself via $g$. This definition extends to the case where $b = 0$ in the $B_{m}$ family by crossing the last root (these are the free-distributions of \cite{meskewnew}). 

There are two more cases of $|2|$-graded $D_{m}$ geometries, given by crossing two of the three extremal roots. Crossing the last two roots results in the free-Lagrangian geometry that is treated in more details in Section \ref{free:la} (though it is formally the same as the previous examples for $b = 1$). The remaining two geometries come from crossing the first root and one of the last two. These geometries are dual to each other, have $G_0 = (\mbb{C}^*)^2 \times SL(m-1,\mbb{C})$ and their real forms are described in section \ref{con:spin} as `conformal spin' geometries.

The second cohomology is slightly simpler than in the $A_m$ family. For most $B_m$, if the last root is crossed, there is one cohomology component, of homogeneity one (there are the `free $m$-distributions'). If any other root is crossed, there are two cohomology components. If the second (contact) root is crossed, both components are of homogeneity one. If the third root is crossed, one component is of homogeneity one and the other of homogeneity zero. In all other cases, they are both of homogeneity zero. 

The two exceptions are $B_2$ and $B_3$ with the last root crossed, which have a single cohomology component, of homogeneity $3$ (see \cite{meskewnew} for more details on the $B_3$ case).

For the $D_m$, $m>5$, conformal-spin geometries have three cohomology components, of homogeneities $1$, $0$ and $0$. Free-Lagrangian geometries, and geometries given by crossing the branching root (the one with three connections) have three cohomology components, all of homogeneity zero. Crossing any other roots results in two cohomological components; if the second (contact) root is crossed, they are both of homogeneity one. If the third root is crossed, one component is of homogeneity one and the other of homogeneity zero. In all other cases, they are both of homogeneity zero.

The algebra $D_5$ nearly follows the above pattern, the only subtlety is that the third root is also the branching root, and has three cohomology components, of homogeneity $1$, $0$ and $0$.

The algebra $D_4$ behaves quite differently; in that case, free-Lagrangian and conformal-spin geometries are isomorphic, and have three components of homogeneity $1$, $1$ and $0$. Crossing the branching (contact) root gives three cohomological components, all of homogeneity one.

\subsubsection{Symplectic geometries}
These are all of the same type, given by $G = Sp(2m, \mbb{C})$ and crossing any root apart from the last one. The $G_0$ is $\mbb{C} \times SL(a,\mbb{C}) \times Sp(2b,\mbb{C})$ with $a+b=m$ and $b \geq 1$.

The geometry is described by two vector bundles $H$ and $K$ of rank $a$ and $2b$ with a skew form $\omega$ on $K$, and
\be
T_{-1} &=& H \otimes K, \\
T_{-2} &=& H \odot H \\
\mc{T} &=& H^* \oplus K \oplus H.
\ee
The algebra\"ic bracket is given by the symmetric tensor product on $H$ and the contraction of $K$ with $\omega$.

The second cohomology is very simple in this case: for $C_m$, $m>2$, the contact (first root crossed) geometries have a single cohomology component, of homogeneity two. All other roots crossed generates two cohomology components. If the second root is crossed, one is of homogeneity $2$, the other of homogeneity $0$. If it's the second-to-last, one component is of homogeneity $1$, the other of homogeneity $0$. In all other cases, they are of homogeneity zero.

The one exception is $C_3$, where the second root is also second-to-last. Here the two components are of homogeneity $2$ and $1$.

\subsubsection{Exceptional geometries} \label{excep:geo}
The exceptional groups have several $|2|$-graded geometries.

\begin{itemize}
\item The group $G_2^{\mbb{C}}$
\end{itemize}

The group $G_2^{\mbb{C}}$ has a single such geometry, given by crossing the long root, with $G_0 = GL(2,\mbb{C})$. This is a contact geometry, given by a vector bundle $H$ of rank two (carrying a natural skew-form $\omega$) and
\be
T_{-1} &=& \odot^3 H \\
T_{-2} &=& \odot^3 (H \wedge H) \\
\mc{T} &=& H \oplus \odot^2 H \oplus H.
\ee
The algebra\"ic bracket is the natural contraction of $\odot^3 H$ with itself via $\omega$. It has a single cohomology component, of homogeneity one.

\begin{itemize}
\item The group $F_4^{\mbb{C}}$
\end{itemize}

The group $F_4^{\mbb{C}}$ carries two such geometries, given by crossing the first or the last root. If the first root is crossed, $G_0 = \mbb{C}^* \times Spin(7,\mbb{C})$, and $T_{-1}$ is the spin representation of $G_0$. Then $T_{-2}$ is seven dimensional, and is the unique seven dimensional irreducible piece in $T_{-1} \wedge T_{-1}$. This is the so-called ``Octonionic contact'' geometry of section \ref{qua:con}. Its Tractor bundle is
\be
\mc{T} = \mbb{C} \oplus T_{-1} \oplus (\mbb{C} \oplus T_{-2}) \oplus T_{-1} \oplus \mbb{C},
\ee
and it has a single cohomology component, of homogeneity zero.

If the last root is crossed, we have a contact geometry with $G_0 = \mbb{C}^* \times Sp(6,\mbb{C})$. The data is given by a vector bundle $H$ of rank $6$ and
\be
T_{-1} &=& \wedge_0^3 H,
\ee
where $ \wedge_0^3 H$ denotes the trace-free subbundle of $\wedge^3 H$. The line bundle $T_{-2}$ is given by the projection $\wedge^2 (\wedge^3_0 H)$ to the line bundle $\wedge^6 H$. Its Tractor bundle is
\be
\mc{T} = H \oplus (\wedge^2_0 H) \oplus H,
\ee
and it has a single cohomology component, of homogeneity one.

\begin{itemize}
\item The group $E_6^{\mbb{C}}$
\end{itemize}

The group $E_6^{\mbb{C}}$ carries four such geometries, though two of them are dual to each other. One is given by crossing the first and the last root, the second by crossing the second or second to last root, and the third by crossing the extremal central root (the one at the end of its own short link).

If we cross the first and the last root, we have $G_0 = (\mbb{C}^*)^2 \times Spin(8,\mbb{C})$. Here the geometry is given by two bundles $H_+$ and $H_-$ of rank $8$, corresponding to the two spin representations of $Spin(8,\mbb{C})$, with
\be
T_{-1} &=& H_+ \oplus H_- \\
T_{-2} &=& H \\
\mc{T} &=& \mbb{C} \oplus H_+ \oplus (H \oplus \mbb{C}) \oplus H_- \oplus \mbb{C},
\ee
where $H$ corresponds to the standard representation of $Spin(8,\mbb{C})$. The algebra\"ic bracket is simply the projection from $H_+ \otimes H_-$ to the single irreducible $H$ component within it. It has three cohomology components, all of homogeneity zero.

If we cross the second root, we have $G_0 = \mbb{C}^* \times SL(2,\mbb{C}) \times SL(5,\mbb{C})$. This is the only geometry that cannot be figured out from simple considerations of the dimension of modules and the decomposition of their symmetric or skew tensor products: there are two valid candidates. However, a detailed look at all the roots of $E_6^{\mbb{C}}$ and the multiplicity of the simple roots inside them resolves the issue: the geometry is given by two bundles $H$ and $K$, of rank $5$ and $2$, such that
\be
T_{-1} &=& (H \wedge H) \otimes K \\
T_{-2} &=& (K \wedge K) \otimes H^* \\
\mc{T} &=& K \oplus (H \wedge H) \oplus (H^* \otimes K) \oplus H
\ee
The algebra\"ic bracket is given by the anti-symmetric volume form of $K$ and the projection onto the irreducible subbundle $\wedge^4 H \cong H^*$ in $(H \wedge H) \odot (H \wedge H)$. It has two cohomology components, of homogeneity zero.

If we cross the extremal root, we have $G_0 = \mbb{C}^* \times SL(6,\mbb{C})$, and the geometry is a contact geometry. It is given by a bundle $H$ of rank $6$ such that
\be
T_{-1} &=& \wedge^3 H \\
T_{-2} &=& \wedge^6 H \\
\mc{T} &=& H^* \oplus \wedge^2 H \oplus H^*.
\ee
The algebra\"ic bracket is the standard map $\wedge^2 (\wedge^3 H) \to \wedge^6 H$. It has a single cohomology component, of homogeneity one.

\begin{itemize}
\item The group $E_7^{\mbb{C}}$
\end{itemize}

The group $E_7^{\mbb{C}}$ carries three such geometries, given by crossing the first root, the second-to-last root, or the extremal root. Crossing the first root gives $G_0 = \mbb{C}^* \times  Spin(12,\mbb{C})$, and the geometry is a contact geometry. The bundle $T_{-1} = H_+$ is given by the rank $32$ spin representation of $G_0$, and there is a one-dimensional module inside $T_{-1} \wedge T_{-1}$, which is isomorphic with $T_{-2}$. If $H_-$ and $H$ are the bundles given by the other spin representation and the standard representation of $Spin(12,\mbb{C})$ respectively, then the Tractor bundle is
\be
\mc{T} = H \oplus H_- \oplus H,
\ee
(note that an outer automorphism of the group interchanges $H_-$ and $H_+$). It has a single cohomology component, of homogeneity one.

Crossing the second to last root gives $G_0 = \mbb{C}^* \times Spin(10,\mbb{C}) \times SL(2,\mbb{C})$. The geometry is defined by two bundles $H_s$ and $K$, the first a spin representation of rank $16$, the second the standard representation of rank $2$, such that
\be
T_{-1} &=& H_s \otimes K \\
T_{-2} &=& (K \wedge K) \otimes H \\
\mc{T} &=& K \oplus H_s \oplus (K \otimes H) \oplus H_s^* \oplus K,
\ee
where $H$ is the rank $10$ bundle coming from the standard representation of $Spin(10,\mbb{C})$. The algebra\"ic bracket is given by skew-symmetrisation on $K$ and the projection from $H_s \odot H_s$ to $H$. It has two cohomology components, both of homogeneity zero.

Crossing the extremal root gives $G_0 = \mbb{C}^* \times  SL(7,\mbb{C})$, and the geometry is defined by a rank $7$ bundle $H$ such that
\be
T_{-1} &=& \wedge^3 H \\
T_{-2} &=& \wedge^6 H \cong H^* \\
\mc{T} &=& H^* \oplus \wedge^2 H \oplus \wedge^2 H^* \oplus H.
\ee
The algebra\"ic bracket is given by the natural map $\wedge^2 (\wedge^3 H) \to \wedge^6 H$. It has a single cohomology component, of homogeneity zero.

\begin{itemize}
\item The group $E_8^{\mbb{C}}$
\end{itemize}

The group $E_8^{\mbb{C}}$ carries two such geometries, given by crossing the first or the last root. For $E_8^{\mbb{C}}$, the standard representation is the adjoint representation; hence the standard Tractor bundle $\mc{T}$ is the adjoint bundle, and for both geometries:
\be
\mc{T} = \ T_{-2} \ \oplus T_{-1} \ \oplus \ \mf{g}_0(T) \ \oplus \ (T_{-1})^* \ \oplus \ (T_{-2})^*.
\ee

If the first root is crossed, $G_0 = \mbb{C}^* \times Spin(14,\mbb{C})$. This gives $T_{-1}$ as the rank $64$ spin representation of $Spin(14,\mbb{C})$. Basic representation theory implies that $T_{-1} \wedge T_{-1}$ contains an irreducible summand of rank $14$, corresponding to the standard representation of $G_0$, which is isomorphic to $T_{-2}$. It has a single cohomology component, of homogeneity zero.

If the last root is crossed, $G_0 = \mbb{C} \times E_7^{\mbb{C}}$. This is a contact geometry, with $T_{-1}$ of rank $56$ as the standard representation of $E_7^{\mbb{C}}$. The group $E_7^{\mbb{C}}$ preserves a skew-form on its standard representation, giving the projection $T_{-1} \wedge T_{-1}$ to the line bundle $T_{-2}$. It has a single cohomology component, of homogeneity one.

\subsection{Non-integrable quaternionic- and octonionic-contact structures} \label{qua:con}
Quaternionic-contact structures are constructed from the $C_m$ families by crossing the second root of the Dynkin diagram:
\be
\textrm{\Cquater}
\ee
The quaternionic-contact aspect is derived from choosing the real form of $G$ to be $Sp(p+1,q+1)$. This results in $G_0$ being $\mbb{R}^* \times Sp(1) \times Sp(p,q)$, while the geometry derives from two left-quaternionic bundles $H$ and $K$, $H$ of quaternionic rank $1$ and $K$ of quaternionic rank $p+q$, with a hermitian metric $h$ of signature $(p,q)$ on $K$. The bundle $H$ has structure bundle $\mc{S}$, an $Sp(1)$-principal bundle, and there is a real rank three vector bundle
\be
\mc{I} = \mc{S} \times_{\rho(\mbb{H})} im \ \mbb{H},
\ee
where $\rho$ denotes the conjugate action of $\mbb{H}$ on $im \ \mbb{H}$. The geometry of these manifolds is given by
\be
T_{-1} &=& \overline{H} \otimes_{\mbb{H}} K \\
T_{-2} &=& \mc{I} \\
\mc{T} &=& H \oplus K \oplus \overline{H}.
\ee
The algebra\"ic bracket is given by the imaginary part of $h$ on $K$ and by the natural contraction of $H$ with itself.

For $m \geq 3$, the second cohomology group of quaternionic contact geometry has two components: one, of homogeneity two, is the obstruction to flatness for an integrable quaternionic-contact structure. The second is the obstruction to integrability for a given distribution. This piece is of homogeneity one for $m = 3$ and of homogeneity zero for $ M \geq 4$.

For this reason, the $m = 3$ case has often been treated differently to others, as being the only case where there existed regular, non-integral distributions $T_{-1}$ generating the quaternionic-contact structure. Our approach allows for the existence of non-integrable $T_{-1}$ distributions in higher dimensions, with the quaternionic structure being in some sense the `best fit' for the distribution.

Similarly, there is an octonionic contact structure, derived from
\be
\textrm{\Focto}
\ee
(see $F_4^{\mbb{C}}$ in section \ref{excep:geo} for more details).

Here, the only component of the cohomology is of homogeneity zero, so the geometry is flat, if regular. Paper \cite{meOlivier} deals with the interesting properties of non-regular quaternionic- and octonionic-contact geometries as conformal infinities of Einstein metrics on the ball.

\subsection{Free-CR and free-Lagrangian geometries} \label{free:la}
These geometries are derived from the $D_m$ family by crossing the last two roots in the Dynkin diagram:
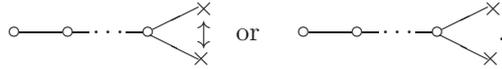
\be
\textrm{\DfreeCR  \ \ \ or \ \ \  \DfreeLan}.
\ee
The above choice of crossed nodes allow three possibilities for $G$: the complex group $G^{\mbb{C}} = SO(2m \mbb{C})$ and the two real forms $SO(m-1,m+1)$ and $SO(m,m)$, corresponding to the two diagrams above. These real forms are entitled free-Cr and free-Lagrangian, respectively. The reason for these names is that for $m = 3$, these are the standard contact-CR and contact-Lagrangian structures, and that in general, the algebra\"ic bracket is free, subject only to the CR or Lagrangian constraints.

The free-CR structure has $G_0 = \mbb{C}^* \times SL(m-1, \mbb{R})$ and is derived from a bundle $E$ of rank $m-1$ and a trivial complex line bundle $L^{\mbb{C}}$ such that
\be
T_{-1} &=& E\otimes L^{\mbb{C}} \\
T_{-2} &=& E \wedge E \\
\mc{T} &=& E^* \oplus (L^{\mbb{C}}) \oplus E.
\ee
If $a$ and $b$ are sections of $L^{\mbb{C}}$, then we may form the symmetric product $Re (a \overline{b})$. This, combined with the identity on $E\wedge E$, gives the algebra\"ic bracket $\wedge^2 T_{-1} \to T_{-2}$.

The free-Lagrangian structure has $G_0 = (\mbb{R}^*)^2 \times SL(m-1,\mbb{R})$ and is derived from a bundle $E$ of rank $m-1$, and a real line bundle $L$, such that
\be
T_{-1} &=& E\otimes L \oplus E \otimes L^* \\
T_{-2} &=& E \wedge E \\
\mc{T} &=& E^* \oplus (L \oplus L^*) \oplus E
\ee
The algebra\"ic bracket derives from the contraction of $L$ with $L^*$, and the identity on $E \wedge E$.

In general there are three cohomology pieces. For $m=4$, two are of homogeneity one and one is of homogeneity zero. For $m>4$, all three are of homogeneity zero.

The free-CR geometry has two pieces of cohomology: the first measures the failure of the bracket to be hermitian with respect to complex structure on $L^{\mbb{C}}$. For $m \geq 3$, it is of homogeneity zero. The second piece measures the failure of spaces of the type $E \otimes l$ to be isotropic, for local sections $l$ of $L^{\mbb{C}}$. It is of homogeneity one for $m=3$ and of homogeneity zero for $m \geq 3$.

The free-Lagrangian has three pieces of cohomology: one is the obstruction to the vanishing of brackets of the type $[e\otimes l, e\otimes l^*]/T{-1}$ for local sections $l$ of $L$ and $e$ of $E$. It is of homogeneity zero for $m \geq 3$. The second measures the failure of $E \otimes L$ to be isotropic, and the third measure the failure of $E \otimes L^*$ to be isotropic. For $m=3$, these two are of homogeneity one, and for $m \geq 4$, they are of homogeneity zero.

The most interesting thing about the free-Lagrangian geometry is that it is a correspondence space for almost spinorial geometries. Almost spinorial geometries are defined by a vector bundle $U$ of rank $m$ with $T = U \wedge U$. They are given by $D_m$ and crossing one of the last two roots. Apart from a few low dimensional exceptions, they have a single piece of harmonic curvature, of homogeneity one.

The correspondence goes as follows: let $N$ be an almost spinorial manifold with $TN = U \wedge U$. Identify the manifold $M$ with the total space of the projectivisation of the bundle $U$. Then $M$ is the correspondence space for $N$.

To see that, let $x$ be a point of $M$. We may identify the vertical vectors of the projection $\pi: M \to N$ with the distribution $E\otimes L^*$, of rank ${m-1}$.

Let $x$ be a point in $M$ and define $u = \pi(x) \in N$. Then $x$ corresponds to a line $R_u$ in $U_{u}$. Wedging this line with $U_{u}$ gives a subspace $K$ of $TN_{u}$, of rank ${m-1}$, and isomorphic with $R_u \otimes (U_{u}/R_u)$. Considerations of the tangent spaces of projective spaces then identifies the vertical component $V_x$ of $T_x$ with the space $R_u^* \otimes (U_{u}/R_u)$. If we define $L_x$ as the pull back of the space $R_u$ and $E_x$ as the pull back of the space $(U_{u}/R_u)$, we get an evident isomorphism $V_x = E_x \otimes L_x$ and $K = D\pi (E_x \otimes L^*_x)$. This defines the bundle $T_{-1}$ at $x$, and the bundle $(T_{-2})_x = T_x / (T_{-1})_x = TN_u / K$ is trivially identified with $\wedge^2 (U_u / R_u) \cong \wedge^2 E_u$.

Of the three pieces of harmonic curvature, two will vanish -- that measuring the isotropy of $E \otimes L^*$ and that obstructing the vanishing of $[e\otimes l, e\otimes l^*]/T{-1}$ -- and all that will remain is the piece measuring the failure of $E$ to be isotropic, which will correspond exactly to the harmonic curvature on $N$ for the almost spinorial structure.

\subsection{Conformal-spin geometries} \label{con:spin}

This geometry may be represented by crossing the first root and one of the last two roots of $D_m$-type algebra. Because of the roots that are crossed, only the complex form $G^{\mbb{C}} = SO(2m,\mbb{C})$ and the fully split form with $G = SO(m,m)$ admits a parabolic of this type. For the split $G$, the $G_0$ is $(\mbb{R}^*)^2 \times SL(m-1)$, and the geometry is given by a single bundle $U$ of rank $m-1$ such that
\be
T_{-1} &=& U^* \oplus \wedge^2 U\\
T_{-2} &=& U \\
\mc{T} &=& \mbb{R} \oplus U^* \oplus U \oplus \mbb{R}.
\ee
The algebra\"ic bracket is the evident contraction between $U^*$ and $\wedge^2 U$. It has three pieces of harmonic curvature, two of homogeneity zero, one of homogeneity one.

What is interesting about this geometry is that it is a correspondence space for both conformal split-signature structures, and almost spinorial structures, as follows:

Let $N$ be a manifold of dimension $2m-2$, and let $[g]$ be a split-signature conformal structure on $N$. Let $Gr(TN)_{m} \to N$ be the Grassmannian bundle of $m-1$-planes in $TN$, and let $M \subset Gr(TN)_{m}$ be the subbundle consisting of those $m-1$-planes that are isotropic with respect to $g \in [g]$. Being isotropic is a conformally invariant concept, so the choice of $g$ does not matter.

For $x \in N$, let $V(t) \in M_x$ be a path for $t \in [0,1]$. Each $V(t)$ is an isotropic $m-1$-plane in $TN_x$. Pick elements $X$ and $Y$ in $V(0)$ and let $\Theta(t)$ be a continuous family of automorphisms of $TN_x$ mapping $V(0)$ to $V(t)$.

Then the isotropy condition implies that
\be
0 &=& g(\Theta(X),\Theta(Y)) \\
&=& g(X,Y) + t g(D\Theta(X)_0,Y) + t g(X,D\Theta(Y)_0) + O(t^2) \\
&=& t \left( g(D\Theta(X)_0,Y) +  g(X,D\Theta(Y)_0) \right) + O(t^2).
\ee
This implies that $D\Theta_0$ is an element of $ \wedge^2 V^*(0)$. Writing $U = V^*$, this we define $(T_{-1})_{V(0)}$ as the subbundle whose projection to $TN$ is $U^*$. Then $T_{V(0)}/(T_{-1})_{V(0)}$ is naturally identified with $TN_x / V(0) = U$, and it is a simple matter to check that the bracket is of the correct form.

For the second correspondence, now let $N$ be an almost-spinorial manifold with bundle $V$ of rank $m$ such that $V \wedge V = TN$. Then we will identify $M$ with the total space of the projectivisation of the bundle $V^*$ (contrast this with the correspondence construction in Section \ref{free:la}). Let $\pi: M \to N$ be the natural projection.

For any point $x$ in $M$, we have a line subspace of $L_x \subset V^*_{\pi(x)}$. This identifies a subspace $L_x^{\perp} \subset V_{\pi(x)}$ of rank $m-1$, which we will call $U$. Basic results on projective spaces then imply that we may identify the vertical tangent space of $TM_x$ with $V^*_{\pi(x)} / L_x$, which is conjugate to $(L_x^{\perp})^* = U^*$. Furthermore $U \wedge U$ is a subspace of $TN_{\pi(x)}$, and we identify $(T_{-1})_x$ with the space that projects down to this space under $D\pi$.

\section{Case Study: $6$-in-$9$}

This section will look into rank $6$ distributions on a $9$-dimensional manifold $M$, demonstrating that partial regularisation will work generically, even in the real category. Furthermore, it will give a computable sufficient condition for when partial regularisations will work.

The parabolic structure used in this situation is the free-Lagrangian structure for $D_4$ (see section \ref{free:la}). By triality, this is also the conformal spin structure for $D_4$ (see section \ref{con:spin}). In terms of crossed Dynkin diagrams, this is
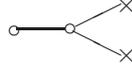
\be
\Dfourcont.
\ee
The real parabolic pair is $(\mf{g}, \mf{p})$, where the graded decomposition of $\mf{g}$ is
\be
\mf{g} = \big( \wedge^2 \mbb{F}^3 \big) \oplus \big( \mbb{F}^3 \oplus \mbb{F}^3 \big) \oplus \big( \mf{gl}(3,\mbb{F}) \oplus \mbb{F} \big) \oplus \big( \mbb{F}^{3*} \oplus \mbb{F}^{3*} \big) \oplus \big( \wedge^2 \mbb{F}^{3*} \big),
\ee
where $\mbb{F} = \mbb{R}$ or $\mbb{C}$. Hence $\mf{g}_0 = \mf{gl}(3,\mbb{F}) \oplus \mbb{F}$ and we may write $\mf{g}_{-1}$ as $E \oplus F$, where both $E$ and $F$ are conjugate to $\mbb{F}^3$. The action of $\mf{sl}(3,\mbb{F}) \subset \mf{g}_0$ gives an isomorphism $E \cong F$ up to scale. For simplicity's sake, pick a scale to fix that isomorphism. Meanwhile, $1 \in \mf{gl}(3,\mbb{F}) \subset \mf{g}_0$ acts by $+1$ on $E$ and $-1$ on $F$, while $1 \in \mbb{F} \in \mf{g}_0$ acts by $+1$ on both these spaces.

Recall (\cite{me2grad1}) that the space of formal brackets $W$ was defined as $\wedge^2 (\mf{g}_{-1})^* \otimes \mf{g}_{-2}$. There is a special element $k$ of $W$, the standard algebra\"ic bracket. To see what it is, pick a basis $\{e_1,e_2,e_3\}$ of $E$ and the corresponding basis $\{f_1,f_2,f_3\}$ of $F$. Let $\{e^1,e^2,e^3\}$ and $\{f^1,f^2,f^3\}$ be the dual bases of $E^*$ and $F^*$.

Since $\mf{g}_{-2} = E \wedge E = F \wedge F$, there is a corresponding basis $\{g_{1},g_{2},g_{3}\}$ of $\mf{g}_{-2}$, given by the relation $g^q = (\epsilon_{ijq}) e_i \wedge e_j$. Here $\epsilon_{ijq}$ is the alternating tensor. Then $k$ itself is defined as:
\be
k &=& \epsilon_{ijq} (e^i \wedge f^j + f^i \wedge e^j) \otimes g_{q}.
\ee
There is a dual bracket $k^*$ in $W^*$, defined as
\be
k^* &=& \epsilon_{ijq} (e_i \wedge f_j + f_i \wedge e_j) \otimes g^{q}.
\ee

Now let $S = GL(\mf{g}_{-1}^*) \times GL(\mf{g}_{-2})$; there is an obvious action of $S$ on $W$. There is a $\partial^*$ operator from $W$ to $\mf{s}$, the Lie algebra of $S$. Define $\alpha$ as the composition of $\partial^*$ with the projection $\mf{s} \to \mf{gl}(\mf{g}_{-1}^*)$ and $\beta$ as the composition of $\partial^*$ with the projection $\mf{s} \to \mf{gl}(\mf{g}_{-2})$. Then $\partial^* = \alpha + \beta$ and if $l$ is an element of $W$, then
\beqa \label{partial:k}
(\alpha l)^i_j &=& \sum_{Iq} l_{qi}^I(k^*)^{qj}_I \\
(\beta l)^I_J &=& -\frac{1}{2} \sum_{ij} l_{ij}^I(k^*)^{ij}_J,
\eeqa
where lower-cases indexes denote $\mf{g}_{-1}$ components while upper-case ones denote $\mf{g}_{-2}$ components.

Define $\theta(l) = \partial^* l - \partial^* k$. Paper \cite{me2grad1} demonstrates that there is an open, $S$-invariant, subset $U \subset W$ such that for all $l \in U$, there exists an $s \in S$ such that $\theta(s \cdot l) = 0$, and the $S$ orbit of $s \cdot l$ is maximally transverse to the kernel of $\theta$ (meaning that locally, $S \cdot (s \cdot l) / \textrm{ker } \theta \cong W/\textrm{ker } \theta$).

Now the bundle $\wedge^2 H^* (\otimes TM/H)$ has structure group $S$, and so can be identified with
\be
\mc{S} \times_S W,
\ee
for the obvious frame bundle $\mc{S}$. Therefore the Levi-bracket $\mc{L}$ defines a function $f_{\mc{L}}$ from $\mc{S}$ to $W$. Paper \cite{me2grad1} demonstrates that whenever the image of $f_{\mc{L}}$ is in $U$ (a well defined property on $M$, since $U$ is $S$-invariant) then there exists a finite choice of normal Cartan connections encoding $H$.

What this section will do is construct $U$, and demonstrate that it is open, dense, and maximum measure in $W$ -- hence that generic $H$ will have the required normal Cartan connections encoding $H$ almost everywhere in $M$. To do this, we will need to decompose $W$ into $S$ orbits -- a task that has already been done, in the complex category, thanks to the work of L.~Yu.~Galitski and D.~A.~Timashev \cite{metabelian}.

\subsection{$S$ orbit classes in $W$}
For this section, assume $\mbb{F} = \mbb{C}$.

The paper \cite{metabelian} actually looks at the class of $S' = SL(\mf{g}_{-1}^*) \times SL(\mf{g}_{-2})$-orbits in $W$. This is not a major difference, however, as $S$ is spanned by $S'$, by scalar multiplication on $W$, and by grading elements of the form $(\lambda^{-2}, \lambda), \lambda \in \mbb{C}^*$. The grading elements fix every element of $W$, so have no effect on orbits. Hence the orbits of $S$ consists of the classes of orbits of $S'$ related by scalar multiplication.

The approach used in \cite{metabelian} is to decompose the lie algebra $\mf{e}_7$ of the complex group $E_7^{\mbb{C}}$ as
\be
\begin{array}{rcccccc}
\mf{e}_7 &=& \mf{h}_{-1} & \oplus &\mf{h}_0 & \oplus & \mf{h}_{1} \\
&=& W^* & \oplus &\mf{s}' & \oplus & W,
\end{array}
\ee
with $\mf{s}'$ the Lie algebra of $S'$. This decomposition has the property that $[\mf{h}_i,\mf{h}_j] \subset \mf{h}_{q}$ where $q = i + j \textrm{ mod } 3$. Note that unlike the parabolic decompositions, the spaces $\mf{h}_{\pm 1}$ are not nilpotent.

Then $\mf{h}_{1}$ is decomposed into semisimple and nilpotent pieces, and the different families are classified according to their semisimple parts. There are seven such families; they are defined in terms of three basic elements $u_1, u_2$ and $u_3$. If $\{a^i\}$, $1\leq i \leq 6$ is a basis for $\mf{g}_{-1}^*$ and $\{b_j\}$, $1 \leq j \leq 3$ a basis for $\mf{g}_{-2}$, then
\be
u_1 &=& a^1 \wedge a^2 \otimes b_1 \ + \ a^3 \wedge a^4 \otimes b_2 \ + \ a^5 \wedge a^6 \otimes b_3, \\
u_2 &=& a^5 \wedge a^4 \otimes b_1 \ + \ a^1 \wedge a^6 \otimes b_2 \ + \ a^3 \wedge a^2 \otimes b_3, \\
u_3 &=& a^3 \wedge a^6 \otimes b_1 \ + \ a^5 \wedge a^2 \otimes b_2 \ + \ a^1 \wedge a^4 \otimes b_3.
\ee

\textbf{Family 1}. This family consists entirely of semisimple pieces, and is full measure in $W$. The canonical form is
\be
l = \lambda_1 u_1 + \lambda_2 u_2 + \lambda_3 u_3,
\ee
where the $\lambda_i$ are constrained by the inequalities
\be
\lambda_i & \neq & 0, \\
(\lambda_i^3 - \lambda_j^3) & \neq & 0, \ i\neq j, \\
\left((\lambda_1^3 + \lambda_2^3 + \lambda_3^3)^3 - (3 \lambda_1^3 \lambda_2^3 \lambda_3^3)^3\right) & \neq & 0.
\ee

\textbf{Family 2}. The canonical form of the semisimple part is
\be
l = \lambda u_2 + \mu u_3,
\ee
subject to the constraints
\be
\lambda, \mu, \lambda^3 \pm \mu^3 \neq 0.
\ee

\textbf{Family 3}. The canonical form of the semisimple part is
\be
l = \lambda u_1 + \mu (u_2 + u_3),
\ee
subject to the constraints
\be
\lambda, \mu, \lambda^3 - \mu^3, \lambda^3 + 8 \mu^3 \neq 0.
\ee

\textbf{Family 4}. The canonical form of the semisimple part is
\be
l = \lambda (u_2 + u_3), \ \ \ \lambda \neq 0.
\ee

\textbf{Family 5}. The canonical form of the semisimple part is
\be
l = \lambda (u_3 - u_2), \ \ \ \lambda \neq 0.
\ee

\textbf{Family 6}. The canonical form of the semisimple part is
\be
l = \lambda u_1, \ \ \ \lambda \neq 0.
\ee

\textbf{Family 7}. This family consists entirely of nilpotent elements.
\\

The main result of this section is that:
\begin{theo} \label{orbit:families:class}
The set $U$ consists of the union of families 1, 2 and 5.
\end{theo}
\begin{proof}
Under the substitution $b_i \to g_i$, and
\be
a^1 \to e^2, \ a^2 \to f^3, \ a^3 \to e^3, \ a^4 \to f^1, \ a^5 \to e^1, \ a^6 \to f^2,
\ee
we get, using equation (\ref{partial:k}),
\be
\partial^* u_1 &=& +\frac{1}{2} \partial^* k, \\
\partial^* u_2 &=& 0, \\
\partial^* u_3 &=& -\frac{1}{2} \partial^* k.
\ee
In this basis, in fact, $k = u_1 - u_3$. Now if we act on $\mf{g}_{-1}$ with the element that sends $a^i$ to $a^{i+1}$ (mod $6$), and act on $\mf{g}_{-2}$ with the element that sends $b_j$ to $- b_{j-1}$ (mod $3$), we permute $u_1$ and $u_2$ while fixing $u_3$. We can similarly permute any other pair of the $u_i$'s.

Consequently, by combining permutations and scalar multiplications, we can put the semisimple parts of all families apart from the seventh into a form such that $\partial^* l = \partial^* k$.

Now assume that $\partial^* l = \partial^* k$. Then let $D_l: \mf{s} \to W$ be the derivative of the action of $S$ on $l$. We know that $D_l = \partial_l$ with
\be
\partial_l (s_{-2}, s_{-1}) (x,y) = l(s_{-1}(x),y) + l(s_{-1}(y),x) - s_{-2}(l(x,y)).
\ee
Now if $\partial^* \partial_l$ is of maximal rank from $\mf{s}$ onto the image of $\partial^*$, then $S$ orbit of $l$ must be maximally transverse to the orbit of $\theta$. A sufficient condition is that $\partial^* \partial_l|_{im \ \partial^*}$ be an isomorphism from ${im \ \partial^*}$ to itself. This happens when the determinant of the map does not vanish. If we set
\be
l = \lambda_1 u_1 + \lambda_2 u_2 + \lambda_3 u_3,
\ee
then
\be
\det \left( \partial^* \partial_l|_{im \ \partial^*} \right) = 0,
\ee
is a polynomial equation in the variables $\lambda_i$. Using Mathematica and Gr\"obner bases methods \cite{grobner}, we can simplify it; the simplification is:
\be
({\lambda_1}^3 -{\lambda_2}^3)({\lambda_1}^3 - {\lambda_3}^3)({\lambda_2}^3 - {\lambda_3}^3) =  0.
\ee
The above equality must be satisfied whenever $\lambda_i = \lambda_j$, $i \neq j$ or whenever two of the $\lambda_i$'s are zero. If all $\lambda_i$ are non-zero, then the restrictions on the first family guarantee that that equation cannot be satisfied. If we assume that $\lambda_1 = 0$, the equation further simplifies to
\be
{\lambda_2} {\lambda_3} ({\lambda_2}^3 - {\lambda_3}^3)  = 0.
\ee
Then the restrictions on the second family imply that this equation cannot be satisfiedeither , while for the fifth family, $\lambda_1 = 0$ and $\lambda_2 = -\lambda_3 \neq 0$ does not solve the equation either.

Hence $\partial^* \partial_l$ is of full rank whenever $l$ is the semisimple part of the first, second or fifth family.

From \cite{metabelian}, we know that if $l = \sigma + \nu$, with $\sigma$ semisimple and $\nu$ nilpotent, then $\sigma$ is in the closure of the $S'$ orbit of $l$. Since $\partial^* \partial_l$ being of maximum rank is an open condition, this implies that for any element $l$ in the families 1, 2 and 5, it can be put into a form (via $S$-action) where $\partial^* l = \partial^* k$ and $\partial^* \partial_l$ is of full rank -- consequently, the $S$ orbit of $l$ is maximally transverse to the kernel of $\theta$.

It now remains to establish that no elements of families 3, 4, 6 or 7 can be in $U$. This is derived from the following lemma:
\begin{lemm}
Let $l$ be in one of the families 3, 4, 6 or 7. Then the stabiliser group $H_l \subset S$ of $l$ contains an element that is not conjugate to any element of $G_0$.
\end{lemm}
\begin{lproof}

The seventh family consists entirely of a finite number of nilpotent elements, which means for each such $l$, there is an element $h \in S'$ that simply scales them by $\lambda$. Consequently, $\lambda^{-1} h$ is in $H_l$; but $\lambda^{-1} h$ is in $S$, and not in $S'$, and is not a grading element. Now $S'$ is preserved under conjugation, and $G_0$ modulo the scaling elements is contained in $S'$. Hence $\lambda^{-1} h$ cannot be conjugate to an element of $G_0$.

For the rest, note that most elements of $G_0$ are not contained in the (conjugation invariant) group $SL(\mf{g}^*_{-1})$. The group $G_0 \cap SL(\mf{g}^*_{-1})$ is one dimensional, and consists of elements that scale $E$ by $\lambda$ and $F$ by $\lambda^{-1}$. Though $E$ and $F$ are not defined up to conjugation, there is a conjugation-invariant result: the fact that the above group consists of diagonalisable maps with three eigenvalues $\lambda$ and three eigenvalues $\lambda^{-1}$.

The proof will proceed by constructing and element $h \in H_l$ that is also an element of $SL(\mf{g}^*_{-1})$ but that is either not diagonalisable, or does not have three eigenvalues $\lambda$ and $\lambda^{-1}$.

For the third family, the nilpotent part can be either zero, or
\be
a^5 \wedge a^3 \otimes b_1 \ + \ a^1 \wedge a^5 \otimes b_2 \ + \ a^3 \wedge a^1 \otimes b_3,
\ee
Now define $h$ as sending $a^2$ to $a^2 + a^1$, $a^4$ to $a^4 + a^3$, $a^6$ to $a^6 + a^5$ and fixing all the other basis elements. This definitely fixes the nilpotent part above, as it fixes $a^1$, $a^3$ and $a^5$. Now $h$ fixes $u_1$, since
\be
h(a^1 \wedge a^2 \otimes b_1) = a^1 \wedge a^2 \otimes b_1 + a^1 \wedge a^1 \otimes b_1 = a^1 \wedge a^2 \otimes b_1,
\ee
and similarly with the other terms in $u_1$. Applying $h$ to $u_2$ and $u_3$, we can see that
\be
h(u_2 + u_3) = u_2 + u_3.
\ee
So $h$ fixes the semisimple part in this family, and hence is in $H_l$. Then, by the argument above, $h$ is not conjugate to any element in $G_0$, as it is not diagonalisable.

The fourth family has six different possible nilpotent parts. All of them are made up of components chosen from among these:
\be
a^1 \wedge a^2 \otimes b_1, \ a^3 \wedge a^4 \otimes b_2, \ a^5 \wedge a^3 \otimes b_1, \ a^1 \wedge a^5 \otimes b_2, \ a^3 \wedge a^1 \otimes b_3.
\ee
Now the previous $h$ will fix all of these components, as well as $u_2 + u_3$. Hence $h \in H_l$ for all $l$ in the fourth family, as before.

The sixth family has fifteen possible nilpotent parts (\cite{metabelian}). Eleven of these parts lack $a^2$ in their expression. Then define $h$ as mapping $a^2$ to $a^2 + a^1$, and fixing all other $a^j$'s; $h$ is not conjugate to any element of $G_0$. This will not affect those eleven nilpotent parts, of course, and will also fix $u_1$. Consequently $h \in H_l$. One of the remaining nilpotent parts is $a^1 \wedge a^3 \otimes b_3 \ + a^2 \wedge a^4 \otimes b_3$. This is stabilised by the map sending $a^5$ to $a^5 + a^6$, and fixing the other basis elements, which also preserves the semisimple part.

The last three remaining nilpotent terms are:
\be
\begin{array}{l}
a^1 \wedge a^4 \otimes b_3 \ + a^1 \wedge a^6 \otimes b_2 \ + a^2 \wedge a^3 \otimes b_3 \ + a^2 \wedge a^5 \otimes b_2 \ + a^3 \wedge a^5 \otimes b_1, \\
a^1 \wedge a^4 \otimes b_3 \ + a^1 \wedge a^5 \otimes b_2 \ + a^2 \wedge a^3 \otimes b_3 \ + a^3 \wedge a^6 \otimes b_1, \\
a^1 \wedge a^3 \otimes b_3 \ + a^1 \wedge a^6 \otimes b_2 \ + a^2 \wedge a^4 \otimes b_3 \ + a^2 \wedge a^5 \otimes b_2.
\end{array}
\ee
The first nilpotent part is stabilised by the map sending $a^2$ to $a^2 - a^1$, $a^4$ to $a^4 + a^3$, $a^6$ to $a^6 + a^4$ and fixing the other basis elements. This also fixes $u_1$, and is evidently not conjugate to an element of $G_0$.

The second nilpotent part is similarly stabilised by the map sending $a^2$ to $a^2 + a^1$, $a^4$ to $a^4 - a^3$ and fixing the rest of the basis elements. The third one is stabilised by the map sending $a^2$ to $a^2 - a^1$, $a^3$ to $a^3 + a^4$, $a^6$ to $a^6 + a^5$ and fixing the rest of the basis elements.
\end{lproof}

Now let $l$ be a member of families 3, 4, 6 or 7, and assume $\partial^* l = \partial^* k$ -- equivalently, $\theta(l) = 0$. We know (\cite{me2grad1}) that for $g \in G_0$, we have $\theta(g\cdot l) = 0$. Now fix a one-parameter subgroup of the stabiliser group of $l$, that is transverse to $G_0$, and let $h$ be an element of this subgroup. Consequently
\be
\theta (g \cdot h \cdot l) = \theta(g \cdot l) = 0.
\ee
The space $W$ is of dimension $6\times 3 \times 5/2 = 45$, the kernel of $\theta$ is of dimension $10$, as is the dimension of $G_0$, and the algebra $\mf{s}$ is of dimension $45$. By the above, there exists a local subset of $S$ around the identity, of dimension at least $11$, that will map $l$ in the kernel of $\theta$. There remains at most $45 - 11 = 34$ degrees of freedom to move $l$ transversely to the kernel of $\theta$, but dimension count implies that this is not enough.

Hence the $S$ orbit of $l$ cannot be maximally transverse to the kernel of $\theta$.
\end{proof}

\subsection{Real brackets}

This section will show that the decomposition of $W$ into orbit types as above continues to be true if we restrict $W$ and $S'$ to the real category. The first steps of the proof will be in the complex category, only shifting to the real category right at the end. Whenever a distance function needs to be calculated, assume the space carries a fixed hermitian metric; the results will be independent of the metric chosen.

For given $u_i$'s we define the element $l(x,y,z)$ to be
\be
x u_1 + y u_2 + z u_3.
\ee

\begin{lemm}
For a given (semisimple) $l(\lambda_1,\lambda_2,\lambda_3)$, we may define a continuous map $\gamma$ from a neighbourhood $N$ of $l(\lambda_1,\lambda_2,\lambda_3) \in \mbb{C}^3$ to $W$ by defining
\be
\gamma(x,y,z) = l(x,y,z).
\ee
Then the closure of $S' \cdot \gamma(N)$ contains a neighbourhood of $l(\lambda_1,\lambda_2,\lambda_3)$ in $W$.
\end{lemm}
\begin{proof}
First note that by dimensional count on the stabilisers, the families 2 to 7 form a set of complex co-dimension at least one in $W$. Consequently the first family is dense in $W$. Moreover, if $q$ is in the first family, then the stabiliser group $H_q \subset S'$ of $q$ consists is the direct product of the finite Well group $G_W$ (see \cite{metabelian}) with the one-dimensional subgroup $H_q \subset S'$, where $h(\mu) \in H_q$ is defined by
\be
h(\mu) (a^i) &=& \lambda a^i \ \ i \in \{1,3,5\}, \\
h(\mu) (a^i) &=& \lambda^{-1} a^i \ \ i \in \{2,4,6\}, \\
h(\mu) (b_j) &=& b_j,
\ee
for all $\mu \in \mbb{C} - \{0\}$ Thus the orbit space of any $l(x,y,z)$ in the first family is isomorphic with $Q = S' / (G_W \times H_q)$. Moreover since $l(x,y,z)$ is semisimple, this space must be closed in $W$.

Now assume that $q =l(\lambda_1,\lambda_2,\lambda_3)$ is in the first family -- which consists entirely of semisimple elements. Then if $N$ is a small \emph{closed} neighbourhood of $(\lambda_1,\lambda_2,\lambda_3)$, we may choose $N$ small enough so that every element of the image of $\gamma$ is in the first family, and that no two elements of $N$ are related by the action of $G_W$.

Then $S' \cdot \gamma(N)$ is isomorphic with $Q \times N$ -- which, by dimensional count, must have non-empty interior in $W$ and also forms a neighbourhood of $q$. Since for any $l(x,y,z) \in \gamma(N)$ the space $Q$ is closed in $W$, $S' \cdot \gamma(N)$ must be closed as well, since $N$ is compact.

Now assume that $q$ is semisimple (hence $q = l(\lambda_1,\lambda_2,\lambda_3)$), but not in the first family. Then pick a small $N$, and let $N' \subset N$ be defined such that $\gamma(N')$ consists of those elements in $\gamma(N)$ that are in the first family. Since the other families have a codimension, $\overline{N'} = N$.

Now consider the space $A = S' \cdot \gamma(N')$. Similarly to before, this must be a space with non-empty interior, isomorphic with $(N'/G_W) \times Q$. Now consider $B = \overline{A} - A$. Obviously $q \in B$; less obviously, we can see that $B$ must have a \emph{real} codimension at least two.

The argument is as follows: if $B$ were of real codimension one, then there would be an element $p$ of the first family in $B$ (since the other families are of \emph{complex} codimension at least one). This $p$ would be isomorphic, via $S$ action, to $l(x,y,z)$ for some $x,y,z$, and the set $S' \cdot \gamma_p (N_p)$ would be isomorphic to $N_p \times Q$ for some closed neighbourhood $N_p$ of $(x,y,z)$ in $\mbb{C}^3$. This would intersect with $A$; an element of the intersection would have an $(x',y',z')$ value coming from $\gamma$, and another $(x'',y'',z'')$ value coming from $\gamma_p$. These values have to be related by an element $h$ of $G_w$. Then replacing $\gamma_p$ with $\gamma_p \circ h$ and $N_p$ with $h^{-1}(N_p)$, we can ensure that $(x',y',z') = (x'',y'',z'')$.

On the overlap $N' \cap N_p$, $\gamma$ and $\gamma_p \circ h$ must be equal. Now define the continuous maps $\pi_1$ as the projection $A \to N'$, and $\pi_2$ as the projection $S'\cdot \gamma_p (N_p)$ to $h^{-1}(N_p)$. These maps are equal on $A \cap S'\cdot \gamma_p (N_p)$, meaning that the value of $h^{-1}(x,y,z)$ must be that of a closure point of $N'$. But this is not possible: if $h^{-1}(x,y,z) \in N'$, then $l(x,y,z) \in A$, contradicting the definition of $B$. But if $h^{-1}(x,y,z) \in (\overline{N'} - N') = N - N'$, then $h^{-1}(x,y,z)$ does not correspond to a member of the first family, contradicting the definition of $p$.

Since $B$ is of real codimension at least two, it cannot be a boundary of $\overline{A}$. Thus $q$ is not a boundary point of $\overline{A}$, and the interior of $\overline{A}$ is a neighbourhood of $q$.
\end{proof}

Now we can define a function $f$ that takes each element $s \cdot l(x,y,z)$ to $(x,y,z) \in \mbb{C}^3$. This can be extended to every element, by considering their semisimple part. It is not, however, well define, as different values of $(x,y,z)$ can correspond to the same orbit. These values are related by the Weyl group $G_W$ (see \cite{metabelian}), a finite subgroup of the Weyl group of $E_7^{\mbb{C}}$, with $1296$ elements. Thus the function
\be
f: W &\to& \mbb{C}^3 / G_W \\
f(s \cdot l(x,y,z)) &=& (x,y,z) / G_W,
\ee
is well defined.

\begin{cor} \label{cor:cont}
The function $f$ is continuous.
\end{cor}
\begin{proof}
Let $q \in W$ be semisimple, i.e.~$q = s \cdot l(\lambda_1,\lambda_2,\lambda_3)$. Then assume there is a sequence
\be
s_n \cdot l(x_n, y_n, z_n),
\ee
where each $l(x_n, y_n, z_n)$ is in the first family, and $(x_n,y_n,z_n)$ does not converge to $(\lambda_1,\lambda_2,\lambda_3)$ (modulo $G_W$). Thus there exists, arbitrarily close to $q$, points of the form $s_n \cdot l(x_n, y_n, z_n)$ with $min_{h \in G_W} ||h \cdot (x_n,y_n,z_n) - (\lambda_1,\lambda_2,\lambda_3) || = \epsilon \neq 0$.

But this contradicts the previous Lemma, since we may take $N'$ to be contained in the ball of centre $(\lambda_1,\lambda_2,\lambda_3)$ and radius $\epsilon / 2 $ -- and know that the closure $S' \cdot \gamma(N)$ is a neighbourhood of $q$, and that there is a neighbourhood of $s_n \cdot l(x_n, y_n, z_n)$ consisting of elements with $f$-values close to $(x_n,y_n,z_n)$.

Now assume that $q$ is not semisimple. As before, there must be a sequence of elements in the first family tending to $q$:
\be
s_n \cdot l(x_n, y_n, z_n).
\ee
Let $s \cdot l(x,y,z)$ be the semisimple part of $q$. Then there is also a sequence $t_m$ such that $t_m \cdot q$ tends to $l(x,y,z)$ as $m \to \infty$.

Now let $\epsilon_m = ||t_m \cdot q - l(x,y,z)||$. Since $t_m$ is non-degenerate, there must also exist a $\delta$ such that $B(q,\delta)$, the ball of centre $q$ and radius $\delta$ gets mapped into $B(t_m \cdot q, \epsilon_m)$, the ball of centre $t_m \cdot q$ and radius $\epsilon_m$.

Then there must exist an $\nu$ such that $s_n \cdot l(x_n,y_n,z_n) \in B(q,\delta)$ for $n > \nu$. We may see $\nu$ as a function of $m$. Now consider the sequence
\be
t_m \cdot s_n \cdot l(x_n,y_n,z_n).
\ee
We let $m$ and $n$ go to infinity according to the following rules: if $n$ is less than the maximum of $\nu(m)$ and $\nu(m+1)$, increment $n$. Otherwise, increment $m$. If we set $\epsilon_0 = \infty$, then
\be
||t_m \cdot s_n \cdot l(x_n,y_n,z_n) - l(x,y,z) || < 2 \epsilon_{m-1},
\ee
so $t_m \cdot s_n \cdot l(x_n,y_n,z_n) \to l(x,y,z)$, implying that $x_n \to x$, $y_n \to y$ and $z_n \to z$ (modulo $G_W$) by the result for semisimple elements.
\end{proof}

Inside each family, there is a maximum element: an element whose stabiliser group is one dimensional, and whose orbit is therefore of maximal dimension. Since $W$ is of dimension $3 (6\times 5)/2 = 45$, and $S'$ is of dimension $3^2 -1 + 6^2 -1 = 43$, the orbit of these maximum elements is of dimension $42$ and co-dimension three.

\begin{lemm} \label{last:lemm}
For any given family, let $q$ be the canonical representative of the element with orbit size 43 under $S'$. Let $l(\lambda_1,\lambda_2,\lambda_3)$ be the semisimple part of $q$. Then there is a continuous map $\gamma$ from a neighbourhood $N$ of $(\lambda_1,\lambda_2,\lambda_3) \in \mbb{C}^3$ to $W$ such that $\gamma(\lambda_1,\lambda_2,\lambda_3) = q$, and such that for all $(x,y,z) \in N$, there exists an element $s_{x,y,z}$ such that the semisimple part of $\gamma(x,y,z)$ is
\be
s_{x,y,z} \cdot l(x,y,z),
\ee
Moreover, if $H_q \subset S'$ is the (one-dimensional) stabiliser group of $q$, then locally the neighbourhood of $q$ is
\be
\Lambda \times N / \mc{I}_q,
\ee
where $\Lambda$ is a neighbourhood of the identity in $S'/H_q$ while $\mc{I}_q$ is a finite subgroup of $S'$.

Allowing for a possible redefinition of $\mc{I}_q$, these results are also true in the real category
\end{lemm}

\begin{lproof}
We will only prove this result for the first three families; the proofs in the other cases are similar, and the results there are less useful.

For the first family, $q = l(\lambda_1,\lambda_2,\lambda_3)$ and the map $\gamma$ is given by the obvious map $\gamma(x,y,z) = l(x,y,z)$. Here $\mc{I}_q  = \{Id\}$, as $\gamma$ is injective.

Locally, the orbit space of $q$ can be identified with $\Lambda \subset S'/H_q$. Let us choose a lift $\phi$ of $\Lambda$ into $S'$. Then by shrinking the size of $N$ and of $\Lambda$ as needed, we can ensure that $\phi(\Lambda)$ acts freely on $\gamma(N)$, giving the result by dimensional count -- the dimension of $\gamma(N)$ is three, while the dimension of $\Lambda$ is $42$ and the dimension of $W$ is $45$.

For the second family, $\lambda_1 = 0$, and the canonical representative is
\be
\lambda_2 u_2 + \lambda_3 u_3 \ +\ a^1 \wedge a^2 \otimes b_1  +  a^3 \wedge a^4 \otimes b_2.
\ee
Let $h(x)$ be the transformation
\be
\begin{array}{rcll}
h(x) (a^j) &=& \frac{1}{x} a^j & j = 1,2,3,4, \\
h(x) (a^j) &=& {x^2} a^j & j = 5,6, \\
h(x) (b_i) &=& \frac{1}{x} b_i & i = 1,2,\\
h(x) (b_3) &=& {x^2} b_3.
\end{array}
\ee
Then $h(x)$ stabilises $u_2$ and $u_3$, and
\be
h(x) u_1 &=& \frac{1}{x^3} a^1 \wedge a^2 \otimes b_1 \ + \ \frac{1}{x^3} a^3 \wedge a^4 \otimes b_2 \ + \ x^6 a^5 \wedge a^6 \otimes b_3,
\ee
Consequently $h(x) x^3 u_1 $ tends to the nilpotent part of $q$ as $x \to 0$, and the map $\gamma$ is given by $h(x) l(x^3,y,z)$. Here $\mc{I}_q$ has three elements, consisting of multiplication by the cube roots of unity. Then since $\gamma$ is injective on $N / \mc{I}_q$, the proof for the neighbourhood description proceeds exactly as in the first family.

For the third family, $\lambda_2 = \lambda_3$ and the canonical representative is
\be
\lambda_1 u_1 + \lambda_2 (u_2 + u_3) \ + a^5 \wedge a^3 \otimes b_1 + a^1 \wedge a^5 \otimes b_2 + a^3 \wedge a^1 \otimes b_3.
\ee
Define $h(x)$ to be the transformation
\be
\begin{array}{rcll}
h(x) (a^j) &=& a^j & j \textrm{ odd,} \\
h(x) (a^j) &=& a^j + \frac{1}{x} a^{j-1} & j \textrm{ even,} \\
h(x) (b_i) &=& b_i.
\end{array}
\ee
Then $h(x) u_1 = u_1$ and
\be
h(x) u_2 &=& u_2 + \frac{1}{x} \left( a^5 \wedge a^3 \otimes b_1 + a^1 \wedge a^5 \otimes b_2 + a^3 \wedge a^1 \otimes b_3 \right) \\
h(x) u_3 &=& u_3 - \frac{1}{x} \left( a^5 \wedge a^3 \otimes b_1 + a^1 \wedge a^5 \otimes b_2 + a^3 \wedge a^1 \otimes b_3 \right).
\ee
Consequently, defining $\gamma$ as
\be
h(y-z) l(x,y,z)
\ee
does the trick. Here $\mc{I}_q = \{Id\}$ just as in the first family, and we are done.

These results did not depend upon the use of any complex numbers (apart from the definition of $\mc{I}_q$ for the second family), and hence are true in the real category.
\end{lproof}

We are now ready for the big result:
\begin{theo}
The classification results of \cite{metabelian} are true in the real category.
\end{theo}
\begin{proof}
Fix real basis $\{a^i\}$ and $\{b_j\}$, and corresponding elements $u_1$, $u_2$ and $u_3$. Define $W'$ to be the set of all elements in $W$ that can be generated from the canonical elements of the various families combined with the action of the (real) group $S'$.

Dimensional considerations imply that elements of the form $s \cdot l(x,y,z)$ constitute an open set in $W$; thus $W'$ has non-empty interior. Then since the third and second families form sets of real co-dimension one, and the other fmailies for sets of higher real co-dimension, there must be an element $r$ of $\overline{W'}$, which is from one of the first three families and has a one-dimensional stabiliser in $S'$.

Then there exists a sequence of real numbers $x_n$, $y_n$ and $z_n$, and a sequence of automorphisms $s_n \in S'$, such that $r$ is the limit of the sequence 
\be
s_n \cdot l(x_n, y_n, z_n).
\ee
Now let $q$ be the canonical representative corresponding to $r$. Therefore there exists a complex map $c \in S'_{\mbb{C}}$ such that $c \cdot r = q$. Now by Corollary \ref{cor:cont}, there are limit points $x_n \to x$, $y_n \to y$ and $z_n \to z$, with the semisimple part of $q$ being $l(x,y,z)$. The $x$, $y$ and $z$ must be real, as $x_n$, $y_n$ and $z_n$ are.

By Lemma \ref{last:lemm}, we know that there must exist elements $t_n \in S'$ such that
\be
t_n \cdot l(x_n, y_n, z_n).
\ee
tends to $q$ as $n \to \infty$. Consequently
\beqa \label{conv:sequ}
(c s_n t_n^{-1})(t_n) \cdot l(x_n, y_n, z_n)
\eeqa
tends to $q$ as $n$ tends to infinity (the limit is just $c \cdot r = q$). Then, by Lemma \ref{last:lemm}, once we are close enough to $q$, there must exist a neighbourhood of the identity $\Lambda$ in $S' / H_q$ and a lift $\phi$ of $\Lambda$ into $S'$ such that
\be
c s_n t_n^{-1} = \phi(\tau_n) h_n /\mc{I}_q,
\ee
where $\tau_n$ is a sequence in $\Lambda$ and $h_n$ is a sequence in $H_q$. Now $\tau_n \to Id$ since the sequence in equation (\ref{conv:sequ}) converges. Hence, modulo $H_q$ and modulo $\mc{I}_q$,
\be
c s_n t_n^{-1} \to Id.
\ee
Since $s_n t_n^{-1}$ is realm this implies that we may choose $c$ to be real as well.

Now the isomorphism of Lemma \ref{last:lemm} from a neighbourhood $N_q$ of $q$ to $N \times \Lambda$ implies that $N_q \subset W'$. By $c$ action, $c(N_q) \subset W'$, and $c(N_q)$ is a neighbourhood of $r$. Consequently neither $q$ nor $r$ can be boundary points of $W'$.

But if $\overline{W'} \neq W$, then the boundary of $W'$ must contain points of the types above, by dimensional considerations. Thus $\overline{W'} = W$, and by reusing Lemma \ref{last:lemm} on individual closure points of $W'$, the whole of $W$ must be in $W'$.

\end{proof}

\subsection{Sufficient conditions and Uniqueness}
We have simplified the problem of determining if $l$ is in $U$ (for both the real case and the complex case) to finding out if $l$ is in one of the families 1, 2 and 5. The first question is whether there is a simple way to determine this, given a specific $l$. The second question is whether there is an easy uniqueness procedure that can be used here, in this much more explicit setting.

\begin{prop}
The second question is easy to solve: different representatives of $G_0$ orbits of an element correspond to different semisimple elements $l(x,y,z)$ with $x-z = 1$. We then simply choose the representative with the lowest $|y|$ value, then, in case of a tie, the lowest $|x|$ value, and then the lowest $|z|$ value, to get a uniqueness procedure that works almost everywhere.
\end{prop}

In principle, the first problem is entirely solvable: include $l$ into $\mf{e}_7$ as described by \cite{metabelian}, and hence into $\mf{gl}(\mf{e}_7)$ by the algebra action -- designate this map from $W$ to $\mf{gl}(\mf{e}_7)$ by $\sigma$. Then the family type of $l$ is entirely determined by the generalised eigenspaces and eigenvalues of $\sigma(l) \in \mf{gl}(\mf{e}_7)$. However, calculating eigenvalues means solving a high degree polynomial, and then using such a result to define generalised eigenspaces means complicated numerical calculations that may compound errors. Ideally, we would want to stick to the one eigenvalue whose value is known: namely $0$. Of course, if $l$ is semisimple, then the generalised $0$ eigenspace for $\sigma(l)$ is precisely its kernel.

\begin{prop} \label{ker:prop}
If $l$ is semisimple, then the kernel $K$ of $\sigma(l)$ splits as $K = K_{-1} \oplus K_{0} \oplus K_1$, where $K_i \subset \mf{h}_i \subset \mf{e}_7$. The space $K_0$ is $d$-dimensional, where $d$ is the dimension of the stabiliser group of $l$ in $S'$, while the dimension of $K_{\pm 1}$ is $d+2$.
\end{prop}
\begin{proof}
Recall the decomposition of the algebra $\mf{e}_7$:
\be
\begin{array}{rcccccc}
\mf{e}_7 &=& \mf{h}_{-1} & \oplus &\mf{h}_0 & \oplus & \mf{h}_{1} \\
&=& W^* & \oplus &\mf{s}' & \oplus & W.
\end{array}
\ee
To define the Lie bracket, we need to fix notation. From now on, assume that $x \wedge y = x \otimes y - x \otimes y$, and similarly for higher wedge products. Then there is a obvious map from $\wedge^2 \mf{g}_{-2}$ to $\mf{g}_{-2}^*$: simply map $x \wedge y$ to $z^*$, where
\be
z^*(z) b_1 \wedge b_2 \wedge b_3 = x \wedge y \wedge z.
\ee
Similarly, there is a map from $\odot^2 (\wedge^2 \mf{g}_{-1}^*)$ to $\mf{g}_1$, sending $x \wedge y $ to $z^*$ where
\be
z^*(z) a^1 \wedge a^2 \wedge a^3 \wedge a^4 \wedge a^5 \wedge a^6 = x \wedge y \wedge z.
\ee
Now the first map is skew, while the second is symmetric, and together they define the bracket $\wedge^2 \mf{h}_1 \to \mf{h}_{-1}$. Minus the dual of this map gives the bracket $\wedge^2 \mf{h}_{-1} \to \mf{h}_1$, while the bracket of $\mf{h} = \mf{s}' = \mf{sl}(\mf{g}_{-1}^*) \oplus \mf{sl}(\mf{g}_{-2})$ on $\mf{h}_{\pm 1}$ is given by the natural action of $\mf{s}'$ on $W$ and $W^*$, and the bracket of $\mf{h}_0$ with itself is the same as $\mf{s}'$.

Assume $c_{\pm 1} \in \mf{h}_{\pm 1}$; then the $\mf{sl}(\mf{g}_{-1}^*)$ component of $[c_1, c_{-1}]$ is the trace-free component of the contraction of $c_1$ and $c_{-1}$ along the $\mf{g}_{-2}$ component and the first $\mf{g}_{-1}^*$ component. The $\mf{sl}(\mf{g}_{-2})$ component of $[c_1, c_{-1}]$ is similarly the trace-free component of the complete contraction $c_1$ and $c_{-1}$ along the $\mf{g}_{-1}^*$ components.

Turning back to the proposition, since $l$ is semisimple, $\mf{e}_7$ decomposes into a direct sum of the eigenspaces of $\sigma(l)$; similarly, it must decompose into a direct sum of the eigenspaces of $\sigma(l)^3$. Now $\sigma(l)^3$ must map each $\mf{h}_i$ to itself, as $l \in \mf{h}_{1}$ and $1 + 1 + 1 = 0$ modulo $3$. Note that $\sigma(l)$ and $\sigma(l)^3$ have the same kernel, which is enough to give the splitting result
\be
K = K_{-1} \oplus K_0 \oplus K_1, \ \ \ K_i \subset \mf{h}_i.
\ee

Let $E^{\lambda}_1 \subset \mf{h}_{1}$ be the eigenspace for $\sigma(l)^3$ with eigenvalue $\lambda \neq 0$. Then we know that
\be
E^{\lambda}_1  = \sigma(l)^3 (E^{\lambda}_1) = \sigma(l)^2 \left( \sigma(l) (E^{\lambda}_1) \right) = \sigma(l) \left( \sigma(l)^2 (E^{\lambda}_1) \right).
\ee
Hence the spaces $\sigma(l) (E^{\lambda}_1) \subset \mf{h}_{-1}$ and $\sigma(l)^2 (E^{\lambda}_1) \subset \mf{h}_{0}$ are of same dimension as $E^{\lambda}_1$. Moreover, the same argument applies for $\sigma(l)(E_1)$ and $\sigma(l)^2 (E_1)$ where
\be
E_1 = \oplus_{\lambda \neq 0} E^{\lambda}_1.
\ee
Defining $E_{-1}$ and $E_0$ as above, we can see that
\be
\mf{h}_i = E_i \oplus K_i.
\ee
Since
\be
\sigma(l)^3 \left( \sigma(l) (E_i) \right) = \sigma(l)^4 (E_i) = \sigma(l) \left( \sigma(l)^3 (E_i) \right) = \sigma(l) (E_i),
\ee
we must have
\be
\sigma(l) (E_i) = E_{i+1},
\ee
modulo $3$ for $i$. \emph{A fortiori} each $E_i$ must have same dimension. Hence the difference in dimensions between $K_0$ and $K_1$ is equal to the difference in dimension between $\mf{h}_0$ and $\mf{h}_1$.

Since $K_0$ is the subalgebra of $\mf{s}'$ that fixes $l$, and since $\mf{h}_0 = \mf{s}'$ is of dimension $43$ while $\mf{h}_1 = W$ and $\mf{h}_{-1} = W^*$ are of dimension $45$, the result follow.

\end{proof}

Fortunately, paper \cite{metabelian} calculates the dimension of the stabiliser of each element of $W$, allowing us to use the above result. Consequently:
\begin{prop} \label{semisimple:condition}
If $l$ is semisimple, then a necessary and sufficient condition for $l$ to be in the first family is that $K_1$ be of dimension $3$. A necessary and sufficient condition for $l$ to be in the fifth family is that $K_1$ be of dimension $11$.
\end{prop}
Calculating the dimensions of kernels is algorithmically tractable, and hence this gives a tractable condition for identifying the first and fifth families in the semisimple cases. For the first family, the $K_1$ is easy to identify: it is simply the span of the $u_i$'s. The second family cannot be identified by the above procedure: $K_1$ is of dimension $5$ in both the second and third family.

If $l$ is not semisimple, we identify the $K_i$'s in terms of its semisimple part. The nilpotent part of $l$ commutes with the semisimple part under $\sigma$, so $\sigma(l)$ must act on $K$ as a nilpotent transformation. The same must be true of $\sigma(l)^3$; as we will be making extensive use of this map, it will be useful to have a simple formula for it, at least for its action on $\mf{h}_{1}$.

Define the map $\phi: \otimes^3 W \to End(W)$ by
\be
& \phi(x_1 \wedge y_1 \otimes z_1, \ x_2 \wedge y_2 \otimes z_2, \ x_3 \wedge y_3 \otimes z_3)\\
&= \\
& x_1 \wedge y_1 (x_2 \wedge y_2 \wedge x_3 \wedge y_3) \otimes \left( z_2 \otimes (z_1 \wedge z_3) - \frac{1}{3} z_1 \otimes(z_2 \wedge z_3) \right) \\
& + \\
& \frac{1}{2} \Big( x_2 \otimes x_1 (y_1 \wedge y_2 \wedge x_3 \wedge y_3) + x_2 \otimes y_1 (x_1 \wedge y_2 \wedge x_3 \wedge y_3) \\
& - y_2 \otimes x_1 (y_1 \wedge x_2 \wedge x_3 \wedge y_3) -y_2 \otimes y_1 (x_1 \wedge x_2 \wedge x_3 \wedge y_3) \\
& -\frac{2}{3} (x_1 \wedge y_1) (x_2 \wedge y_2 \wedge x_3 \wedge y_3) \Big) \otimes z_1(z_2 \wedge z_3),
\ee
using the identifications $\wedge^4 \mf{g}_{-1}^* \cong \wedge^2 \mf{g}_{-1}$ and $\wedge^2 \mf{g}_{-2} \cong \mf{g}_{-2}^*$ of Proposition \ref{ker:prop}. Then $\sigma(l)^3|_W$ is $\phi(l,l,l)$.

The use of this construction will be obvious in the following theorem:
\begin{theo}
The kernel of $\sigma(l)^3|_W$ is three dimensional, if and only if $l$ is in the first family.
\end{theo}
\begin{proof}
The `if' part of the statement is trivially true. The `only if' is true for $l$ semisimple, by Proposition \ref{semisimple:condition}. So assume now that $l$ is not in the first family, not semisimple, and that $\sigma(l)^3|_{\mf{h}_{1}}$ has a kernel of dimension equal to three. Note that the last condition is equivalent with $\sigma(l)^3|_{K_1}$ having a kernel of dimension equal to three.

Define $L$ to be the kernel of $\sigma(l)$. Since $l$ is of degree $+1$, $L$ must split as $L_{-1} \oplus L_0 \oplus L_{1}$.
\begin{lemm}
The spaces $L_{\pm 1}$ are of dimension at least three. The space $L_0$ is of dimension at least one.
\end{lemm}
\begin{lproof}
Let $l_n$ be a sequence of elements in the first family tending to $l$. The kernels of $\sigma(l_n)|_{\mf{h}_i}$ are of dimensions $3$, $1$ and $3$ in degrees $-1$, $0$ and $1$ (see \cite{metabelian} or the proof of Proposition \ref{ker:prop}. Then we just use the result from linear algebra that the rank of the kernel of a sequence of linear transformations is upper-semi continuous.
\end{lproof}
For $\sigma(l)^3$ to had a kernel of dimension three, the above dimensional bounds must be sharp, and we must have $\sigma(l)(K_1) \cap L_{-1} = 0$ and $\sigma(l)^2(K_1) \cap L_0 = 0$. Since $\sigma(l)$ is nilpotent on $K_1$, we must have $\sigma(l)^3(K_1) \cap L_{1} \neq 0$.

Since $\sigma(l)(K_i) \subset K_{i+1}$, we know that the space $K_{-1}$ and $K_0$ split as
\be
K_{-1} &=& \sigma(l)(K_1) \oplus L_{-1} \\
K_0 &=& \sigma(l)(K_0) \oplus L_0.
\ee
Now applying $\sigma(l)^3$ to $K_{-1}$, we get the image as $\sigma(l)^4(K_1)$. Since $\sigma(l)^3(K_1) \cap L_{1} \neq 0$, this space must have one dimension less than $K_1$; i.e.~$\sigma(l)^3$ has a kernel of dimension at least four on $\mf{h}_{-1}$.

We now note that $\mf{e}_7$ carries a killing form $B$. It is easy to see that under $B$, $\mf{h}_{\pm 1}$ are isotropic and dual to each other. The form $B$ must be $\mf{h}_0 = \mf{s}'$ invariant: hence it must be the natural contraction between $W$ and $W'$, multiplied by some constant $\lambda \neq 0$.

Since $\sigma(l)$ derives from the algebra action, if must preserve $B$. Consequently if $q^* \in \mf{h}_{-1}$ and $q \in \mf{h}_1$:
\be
q^* \llcorner \sigma(l)^3(q) &=& \frac{1}{\lambda} B(q^*, \sigma(l)^3 q) \\
&=& -\frac{1}{\lambda} B(\sigma(l) q^*, \sigma(l)^2 q) \\
&=& -\frac{1}{\lambda} B(\sigma(l)^3 q^*, q) = -\sigma(l)^3 (q^*) \llcorner q.
\ee
Thus the action of $\sigma(l)^3 $ on $\mf{h}_{-1}$ is minus the dual of its action on $\mf{h}_1$. This, however, is a contradiction, as we have shown that the two have kernels of different dimensions.

Consequently our original assumption was wrong, and there are no $l$'s not in the first family such that $\sigma(l)^3|_{\mf{h}_{1}}$ has a kernel of dimension equal to three.

\end{proof}

\bibliographystyle{amsalpha}
\bibliography{ref}

\providecommand{\bysame}{\leavevmode\hbox to3em{\hrulefill}\thinspace}
\providecommand{\MR}{\relax\ifhmode\unskip\space\fi MR }
\providecommand{\MRhref}[2]{%
  \href{http://www.ams.org/mathscinet-getitem?mr=#1}{#2}
}
\providecommand{\href}[2]{#2}
\begin{thebibliography}{Arm07}

\bibitem[AB]{meOlivier}
Stuart Armstrong and Olivier Biquard, \emph{Einstein metrics with anisotropic
  boundary behaviour}, arXiv, arXiv:0901.1051v1 [math.DG].

\bibitem[Arma]{me2grad1}
Stuart Armstrong, \emph{Non-regular $|2|$-graded geometries {I}: general
  theory}, arXiv, arXiv:0902.1133v1 [math.DG].

\bibitem[Armb]{meCR}
\bysame, \emph{Reducing almost {L}agrangian structures and almost {CR}
  geometries to partially integrable structures}, arXiv, arXiv:0807.1234v2
  [math.DG].

\bibitem[Arm07]{meskewnew}
\bysame, \emph{Free $3$-distributions: holonomy, {F}efferman constructions and
  dual distributions}, arXiv (2007), arXiv:0708.3027v3 [math.DG].

\bibitem[Buc]{grobner}
Bruno Buchberger, \emph{An algorithmic criterion for the solvability of a
  system of algebraic equations}, Aequationes Mathematicae \textbf{4},
  374--383.

\bibitem[{\v{C}}ap05]{capauto}
Andreas {\v{C}}ap, \emph{Automorphism groups of parabolic geometries},
  {`Proceedings of the 24th Winter School on Geometry and Physics, Srni 2004'}
  Rend. Circ. Mat. Palermo Suppl. ser. II \textbf{75} (2005), 233--239.

\bibitem[{\v{C}}G00]{TBIPG}
Andreas {\v{C}}ap and Rod Gover, \emph{Tractor bundles for irreducible
  parabolic geometries}, S.M.F. Colloques, Seminaires \& Congres \textbf{4}
  (2000), 129--154.

\bibitem[{\v{C}}G02]{TCPG}
\bysame, \emph{Tractor calculi for parabolic geometries}, Trans. Amer. Math.
  Soc. \textbf{354} (2002), no.~4, 1511--1548.

\bibitem[{\v{C}}N08]{capkatauto}
Andreas {\v{C}}ap and Katharina Neusser, \emph{On automorphism groups of some
  types of generic distributions}, arXiv (2008), arXiv:0807.0974v1 [math.DG].

\bibitem[{\v{C}}Sed]{capslo}
Andreas {\v{C}}ap and Jan Slov{\'a}k, \emph{{P}arabolic {G}eometries {I}:
  {B}ackground \& general theory}, To be published.

\bibitem[GT99]{metabelian}
L.~Yu. Galitski and D.~A. Timashev, \emph{On classification of metabelian {L}ie
  algebras}, Journal of Lie Theory \textbf{9} (1999), 125--156.

\bibitem[Kos61]{Kostant}
Bertram Kostant, \emph{Lie algebra cohomology and the generalized
  {B}orel-{W}eil theorem}, Ann. of Math. (2) \textbf{74} (1961), 329--387.

\bibitem[Miy91]{liecontact}
Reiko Miyaoka, \emph{Lie contact structures and normal {C}artan connections},
  Kodai Math. J. \textbf{14} (1991), no.~1, 13--41.

\bibitem[Zad]{liecontact2}
Vojtech Zadnik, \emph{Lie contact structures and chains}, arXiv,
  arXiv:0901.4433v1 [math.DG].

\end{thebibliography}

\end{document}